\documentclass{amsart}

\usepackage{amsmath}
\usepackage{paralist}
\usepackage{amsfonts}
\usepackage{MnSymbol}
\usepackage{graphics} 
\usepackage{epsfig} 
\usepackage{psfrag}
\usepackage{pdfsync}
\usepackage[usenames]{color}
\usepackage{cancel}
\usepackage{mathrsfs}
\DeclareMathAlphabet{\mathpzc}{OT1}{pzc}{m}{it}

\definecolor{arancio}{rgb}{0.90,0.50,0.20}

\definecolor{blu}{rgb}{0.,0.,1.}

\definecolor{pavone}{rgb}{0.00,0.00,0.63}

\definecolor{malva}{rgb}{0.10,0.50,0.50}

\definecolor{rosso}{rgb}{1.,0.,0.}

\definecolor{geranio}{rgb}{0.90,0.00,0.20}

\definecolor{cerulean}{rgb}
{0.0, 0.48, 0.65}

\newtheorem{theorem}{Theorem}[section]

\newtheorem{lem}[theorem]{Lemma}
\newtheorem{prop}[theorem]{Proposition}

\theoremstyle{definition}
\newtheorem{definition}{Definition}[section]
\newtheorem{remark}{Remark}[section]

\newcommand{\N}{\mathbb{N}}
\newcommand{\R}{\mathbb{R}}
\newcommand{\re}{\mathbb{R}}

\date{\today}

\newcommand{\bcl}{\begin{center}}
\newcommand{\ecl}{\end{center}}
\newcommand{\brl}{\begin{right}}
\newcommand{\erl}{\end{right}}
\newcommand{\ben}{\begin{enumerate}}
\newcommand{\barr}{\begin{array}}
\newcommand{\earr}{\end{array}}
\newcommand{\btab}{\begin{tabular}}
\newcommand{\etab}{\end{tabular}}
\newcommand{\bdoc}{\begin{document}}
\newcommand{\edoc}{\end{document}}
\newcommand{\beqy}{\begin{eqnarray}}

\newcommand{\beq}{\begin{equation}}
\newcommand{\beqi}{\begin{eqnarray*}}
\newcommand{\bitem}{\begin{itemize}}
\newcommand{\brem}{\begin{remark}}
\newcommand{\erem}{\end{remark}}
\newcommand{\eitem}{\end{itemize}}
\newcommand{\nln}{\newline}
\newcommand{\newt}{\newtheorem}

\renewcommand{\a }{\alpha }
\renewcommand{\b }{\beta }
\newcommand{\g }{\gamma}
\newcommand{\G }{\Gamma }
\renewcommand{\d }{\delta }
\newcommand{\D }{\Delta }
\newcommand{\e }{\epsilon }
\newcommand{\z }{\zeta }
\renewcommand{\l }{\lambda  }
\renewcommand{\L }{\Lambda }
\newcommand{\m }{\mu }
\newcommand{\n }{\nu }
\renewcommand{\r }{\rho }
\newcommand{\s }{\sigma }
\newcommand{\Sig }{\Sigma }
\renewcommand{\t }{\tau }
\newcommand{\vT}{\Theta }
\renewcommand{\o }{\omega }
\renewcommand{\O }{\Omega }
\newcommand{\pa }{\partial  }

\newcommand{\supp}{\text{\rm supp}\,}
\newcommand{\sgn}{\text{\rm sgn}\,}

\newcommand{\red}[1]{{\color{red}{#1}}}
\newcommand{\blue}[1]{{\color{blue}{#1}}}

\title[Extinction and propagation on metric trees] {Extinction and propagation phenomena \\ for semilinear parabolic equations on metric trees}

\author[Punzo]{Fabio Punzo}
\address{Dipartimento di Matematica, Politecnico di Milano \\ Via Bonardi 9, I-20133 Milano, Italy}
\email{fabio.punzo@polimi.it}
\author[Tesei]{Alberto Tesei}
\address{Dipartimento di Matematica ``G. Castelnuovo", Universit\`a Sapienza di Roma\\ P.le A. Moro 5, I-00185 Roma, Italy, and
Istituto per le Applicazioni del Calcolo ``M. Picone", CNR,
Via dei Taurini 19, I-00185 Roma, Italy}
\email{alberto.tesei@lincei.it}



\date{\today}

\begin{document}

\bibliographystyle{h-elsevier2}

\begin{abstract}
We study the Cauchy-Neumann problem on a regular metric tree $\mathcal{T}$ for the semilinear heat equation with forcing term  of KPP type. Propagation and extinction of solutions, as well as asymptotical speed of propagation are investigated. With respect to the corresponding problem in the Euclidean space new phenomena arise, which depend both  on the properties of the diﬀusion process and on the geometry of $\mathcal{T}$. Significant parallels emerge with the analogous problem in the hyperbolic space.
\end{abstract}

\maketitle


\section{introduction}

\subsection{Motivations} 

Metric graphs appear in several areas of science with different motivations (e.g., see \cite{BF, BK, C, MSV, N, SFPP} and references therein). In particular, they provide the proper framework to model spatial structures of interest in physics and biology, like e.g. river basins, cave systems and vascular networks \cite{Ra, SCA, SMA}.

In this framework, dynamics is ruled by partial differential equations stated on the edges of the graph, subject to boundary conditions at the vertices. Therefore, a natural question is how  geometric properties (i.e., features like branching degree of vertices, length of edges, or symmetries) influence the dynamics - an issue already addressed for Riemannian manifolds (e.g., see \cite{BPoT,BPuT, MPT}). In this connection, observe that metric graphs can be regarded as one-dimensional manifolds with singularities.

Compared to the case of manifolds, for metric graphs information concerning spectral theory, or important features like curvature and volume doubling, is less available than in the case of manifolds (e.g., see \cite{PT} and references therein). Therefore, it is arguable that the above question can be fully answered only for specific classes of graphs.

Metric trees (see Definition \ref{detree}) form an interesting class of general metric graphs. Among them, {\em regular trees} have a special metric and combinatorial structure, which is reflected by the geometry of the corresponding $L^2$-spaces and Sobolev spaces. Roughly speaking, due to this structure a problem on a regular tree can be reduced to a family of more elementary problems on intervals; in particular, this gives quite a lot of information concerning  spectral theory (see Subsection \ref{metre}). Both by these features and by their specific modelling interest (e.g., see \cite{SCA, Steta, VCNA}), regular metric trees seem a promising candidate for a case study of the stated problem.


\subsection{Statement of the problem and assumptions}

 In the light of the above remarks, in this paper we study the Cauchy-Neumann problem
 \begin{subequations}\label{CN}
\begin{equation}\label{diffeq}
u_t  = \Delta u+ f(u) \quad \text{in $\mathcal{T}\times \R_+$\,,}
\end{equation}
\begin{equation}\label{inco}
u \, = u_0 \quad \text{in $\mathcal{T}\times \{0\}$}
\end{equation}
\end{subequations}
on a regular metric  tree $\mathcal{T}$, with homogeneous Neumann condition at the root $O$. Here $\Delta\equiv\D_N$ denotes the {\em Neumann Laplacian} on $\mathcal{T}$ (see Subsection \ref{megra}).

The {\em generating sequences} $\{b_n\}$ and $\{\r_n\}$ are important to describe the structure of a regular tree: the integer $b_n\ge2$ denotes the number of edges emanating from any vertex of the $n$-th generation, and $\r_{n+1}-\r_n>0$ is the common length of such edges (see Subsection \ref{metre}). We always make the following assumption: 
$$
\left\{
\begin{array}{l}
(i)\;\;\;
\text{$\mathcal{T}$ is is a regular metric tree {\em of infinite height,}}\\
\qquad\qquad\qquad h(\mathcal{T})\,:=\,\sup_{x\in\mathcal{T}}\r(x)\,=\,\infty\,;\smallskip\\
(ii) \;\;\text{$\{b_n\}$ is  nondecreasing, $\{\r_n-\r_{n-1}\}$ nonincreasing}
\end{array}
\right.
\leqno(H_0)
$$
(here $\r(x):=d(x,O)$ denotes the distance from the root of any $x\in \mathcal{T}$). In particular, $(H_0)$-$(ii)$ holds for {\em homogeneous trees} with $b_n=b,\,\r_n-\r_{n-1}=r$ ($n\in\N$; see Definition \ref{decre}-$(iv)$).

\smallskip

Concerning the Cauchy data $u_0$, the following assumption is made:
$$
u_0\in L^2(\mathcal{T})\,,\qquad 0\le u_0\le 1\;\;\textrm{a.e. in $\mathcal{T}$}.
\leqno(H_1)
$$
By assumption $(H_1)$, $u_0$ also belongs to $L^\infty( \mathcal{T})$ and is of the form
$u_0= \bigoplus_{n=1}^\infty\,\bigoplus_{i=1}^{b_n}\,u_{0i}\,. $

Concerning the function $f$, we always assume that
$$
 f\in C^1([0,1])\,\cap \,C^\infty(0,1)\,,\quad f(0)=f(1)=0\,,\quad f(u)>0\;\;\textrm{for any}\;\; u\in (0,1)\,. \leqno (H_2)
$$
Moreover, the following assumptions will be used:
\begin{equation*}
\sup_{u\in(0,1]}\frac{f(u)}{u}=: M<\infty\,,
\leqno{(H_3)}
\end{equation*}
\begin{equation*}
\sup_{u\in(0,1]}\frac{f(u)}{u}=f'(0)\,,
\leqno{(H_4)}
\end{equation*}
\begin{equation*}
 \limsup_{u\to 0^+}\, u^{-p}f(u)<\infty\quad\textrm{for some\,}\; p >1\,.
\leqno{(H_5)}
\end{equation*}

The function $f$ is said to be {\em of KPP type}, if it satisfies $(H_2)$ and there holds $f'(0)\in \R_+$. Clearly, any function $f$ satisfying $(H_2)$ and $(H_4)$ is of KPP type. Moreover, if $(H_2)$ holds, assumption $(H_4)$ implies $(H_3)$.

\smallskip

Solutions of problem \eqref{CN} are  meant in the sense of Definition \ref{deso} below. By Theorem \ref{wpp} a unique solution $u$ of problem \eqref{CN} exists, if $(H_1)$-$(H_2)$ hold. By definition, for any $t\in\R_+$\,:
\begin{itemize}
\item[-] $u(\cdot,t)$ belongs to the domain of the Laplacian, thus in particular is continuous at any vertex of $\mathcal{T}$;
\item[-] $u(\cdot,t)$ satisfies the {\em Kirchhoff condition} at any vertex of $\mathcal{T}$. In particular, there holds $\frac{\pa u}{\pa \n}(O,t)= 0$ for all $t\in\R_+$, \,$\frac{\pa }{\pa \n}$ denoting the outer normal derivative (see Subsection \ref{lamegra}).
\end{itemize}

Let us mention that the analogue of problem \eqref{CN} on discrete trees, with $f$ of  KPP type,  was considered in \cite{hh}.


\subsection{The Euclidean case}\label{oueu}

 Problem \eqref{CN} can be regarded as the counterpart of the Cauchy problem in $\R^N$
\begin{equation}\label{e1i}
\left\{
\begin{array}{ll}
 u_t  = \Delta u+ f(u)  &\textrm{in}\,\,\R^N\times \R_+
\\& \\
u \, = u_0& \textrm{in}\,\, \R^N\times \{0\}\,,
\end{array}
\right.
\end{equation}
which has been widely investigated (in particular, see \cite{AW2,AW, BH,FML,Jon}). Let us recall some well-known results (see \cite{AW} for details).

\smallskip

\noindent $(a)$  Let the forcing term $f$ satisfy $(H_2)$ and the following condition:
 \begin{equation}\label{exticondir}
 \limsup_{u\to 0^+}\, u^{-p}f(u)<\infty \quad\textrm{for some\,}\; p >1+\frac{2}{N}\,.
\end{equation}
Then, if $u_0$ is small enough, extinction occurs:
 \begin{equation}\label{exti}
\lim_{t\to \infty} u(x,t)\,=\,0\quad\textrm{uniformly in\,}\; \re^N.
\end{equation}
Observe that $\displaystyle{p=1+\frac{2}{N}}$ is the {\em Fujita exponent} of problem \eqref{e1i} with $f(u)=u^p$ (see \cite{Fu66}).

\smallskip

\noindent $(b)$
If $f\in C^1([0,1])$, $f(0)=f(1)=0$, and
$$
\left\{
\begin{array}{l}
(i)\;\;\hbox{there exists}\;a\in (0,1]\;\hbox{such that}\;\, f(u)>0\;\hbox{for any}\;u\in (0,a)\,,\smallskip \\
(ii)\,\displaystyle{\liminf_{u\to 0^+}\, u^{-(1+\frac{2}{N})}f(u)>0}\,,
\end{array}
\right.
\leqno(HT)
$$
then for every solution $u \not \equiv 0$ of problem \eqref{e1i} there holds
\begin{equation}\label{propp}
\liminf_{t\to \infty} u(x,t)\,\ge\, a\quad\textrm{uniformly on compact subsets of\,}\; \re^N\,.
\end{equation}
This is the so-called {\em hair-trigger effect}. In particular, assumption $(HT)$-$(i)$ is satisfied with $a=1$ if $f$ is of  KPP type, thus propagation always occurs in this case - namely, for every solution $u \not \equiv 0$ of problem \eqref{e1i} there holds
\begin{equation}\label{e2i}
\lim_{t\to \infty} u(x,t)\,=\,1\quad\textrm{uniformly on compact subsets of\,}\; \re^N.
\end{equation}
More generally, assumption $(HT)$ is satisfied  if $(HT)$-$(i)$ holds and
\begin{equation*}
 \liminf_{u\to 0^+}\, u^{-p}f(u)>0\quad\textrm{for some\,}\; p \in \left(1,1+\frac{2}{N}\right)\,,
\end{equation*}
thus \eqref{propp} follows also in this case.
\smallskip

\noindent $(c)$ For KPP there exists an {\em asymptotic speed of propagation} $c_0>0$, which is uniquely determined by the following properties:

\noindent $(i)$ no solution of problem \eqref{e1i} with $u_0$ having compact support propagates with speed greater than $c_0$. In fact, for any $c>c_0$ and $y\in\re^N$
\begin{equation}\label{propa1}
\lim_{t\to \infty} \sup_{|x-y|>c t}\, u(x,t)\,=\,0 \,;
\end{equation}

\noindent $(ii)$ if a solution of problem \eqref{e1i} propagates, then its speed is no smaller than $c_0$. In fact, if
\begin{equation*}
\liminf_{t\to \infty} u(x,t)\,\ge\, a\quad\textrm{uniformly on compact subsets of $\R^N$ for some $a \in (0,1]$,}\,
\end{equation*}
 then for any $c <c_0$ and $y\in\re^N$
\begin{equation}\label{propa2}
\liminf_{t\to \infty} \inf_{|x-y|<c t}\, u(x,t)\,\ge\,a \,.
\end{equation}

Remarkably, the asymptotic speed of propagation $c_0$ only depends on the forcing term $f$. In fact, its definition relies on the properties of solutions of the ordinary differential equation
\begin{equation}\label{e84}
q''+c q' + f(q)=\,0\quad \hbox{in } \re \qquad (c \in\re)
\end{equation}
(see \cite[Theorem 4.1 and Lemma 4.3]{AW}). Equation \eqref{e84} arises when seeking {\em plane wave solutions} of the  equation
\begin{equation}\label{e92e}
u_t
 = \Delta u+ f(u) \quad
 \textrm{in}\,\,\re^N\times \R_+
 \end{equation}
 - namely,  solutions of equation \eqref{e92e} of the form
\begin{equation}\label{e7i}
u(x,t)=q(\langle x, \nu\rangle - c t) \qquad (x\in \re^N, t>0) \,,
\end{equation}
where $q$ is a real function,  $\nu\in\re^N$ is a fixed unit direction, $c\in \re$ and $\displaystyle{\langle x, \nu\rangle:= \sum_{i=1}^n\, x_i{\nu}_i\,}$.


\subsection{The hyperbolic case}

Concerning the companion problem of \eqref{CN} in the hyperbolic space $\mathbb{H}^N$ $(N\ge 2)$,
\begin{equation}\label{e1}
\left\{
\begin{array}{ll}
  u_t = \Delta_H u + f(u)  &\textrm{in}\,\,\mathbb{H}^N\times \re_+
\\& \\
\textrm{ }u \, = u_0& \textrm{in}\,\, \mathbb{H}^N\times \{0\}
\end{array}
\right.
\end{equation}
(where $\Delta_H$ denotes the Laplace-Beltrami operator in $\mathbb{H}^N$), a novel feature is the existence of a {\em threshold effect} for extinction prevailing over propagation, which has no counterpart in $\R^N$ where propagation always takes place.  This depends on the size of $f'(0)$ and  is related to the fact that the $L^2$-spectrum of the operator $-\Delta_H$
in $\mathbb{H}^N$ has positive infimum $\l_1:=\frac{(N-1)^2}{4}$\,. In fact, let $u_0$ be continuous in $\mathbb{H}^N$, $0\le u_0(x)\le 1$ for any $x \in \mathbb{H}^N$, and let $(H_2)$ and $(H_4)$ hold. Then the following holds (see \cite[Theorem 3.2]{MPT}):

\smallskip

\noindent - if $f'(0)>\l_1$ and $u_0\not\equiv0$, the solution of the Cauchy problem \eqref{e1}
converges to $1$ uniformly on compact subsets of $\mathbb{H}^N$ as $t\to\infty$;

\smallskip

\noindent - if $f'(0)<\l_1$ and $u_0$ has compact support, as $t\to\infty$ the solution of \eqref{e1} converges to $0$ uniformly in $\mathbb{H}^N$.

 In this connection, observe that if $(H_4)$ holds (thus $f$ is of  KPP type), by \cite[Proposition 4.2]{AW} there holds $c_0=2\sqrt{f'(0)}$, thus
\begin{equation}\label{e303}
 f'(0)\,>\,\l_1\,=\,\frac{(N-1)^2}{4} \qquad \Longleftrightarrow \qquad  c_0\,>\,N-1 \,.
\end{equation}
Therefore, in this case extinction prevails if  $c_0<N-1$ (and $u_0$ has compact support), whereas in the opposite case  the asymptotic speed of propagation is $c_0- (N-1)=2\big[\sqrt{f'(0)} -\sqrt{\l_1}\big]$ instead of $c_0$ (see \cite[Theorem 3.6]{MPT}). Namely, the speed of propagation is less than that in $\R^N$ and depends on the space dimension beside the source term $f$.

In colloquial terms, the positivity of $\l_1$ makes diffusion stronger in $\mathbb{H}^N$ than in $\R^N$. As a consequence, in $\mathbb{H}^N$ the Fujita exponent  is $p=1$, i.e,  the Fujita phenomenon does not occur (in this connection, observe that the volume  growth in $\mathbb{H}^N$ is exponential; see \cite{BPoT}, \cite[Proposition 8.3]{PT}). The enhanced effect of diffusion explains why   in $\mathbb{H}^N$ extinction can prevail against propagation, and
why the asymptotic speed is less than that in $\R^N$ if propagation occurs.


\subsection{Outline of results }

As for extinction versus propagation, the situation for problem \eqref{CN} is strongly reminiscent of that concerning $\mathbb{H}^N$. In fact, under assumptions $(H_1)$, $(H_2)$ and $(H_4)$ the following will be proven:

\noindent - if $f'(0)>E_0$, all solutions of problem \eqref{CN} with $u_0\not\equiv0$ converge to $1$ uniformly on compact subsets of $\mathcal{T}$ as $t\to\infty$ (see Theorem \ref{t3});

\noindent - if $f'(0)<E_0$ (and $\mathcal{T}$ is homogeneous), all solutions of \eqref{CN} with $u_0\not\equiv0$ suitably small converge to $0$ uniformly in $\mathcal{T}$ as $t\to\infty$ (see Theorem \ref{t30}).

\noindent The same conclusion of Theorem \ref{t30} holds - namely, extinction occurs - if assumption $(H_4)$ is replaced by $(H_5)$, which should be compared with condition \eqref{exticondir} of the Euclidean case (see Theorem \ref{t2}).

\smallskip

To explain the outlined analogy observe that, if $(H_1)$ holds,

\noindent $(i)$ the infimum $E_0$ of the $L^2$-spectrum of $-\D$ is strictly positive (see Theorem \ref{tre00} and the spectral analysis in \cite{NS1,NS2,SS,So1,So2});

\noindent $(ii)$ the Fujita phenomenon does not occur in regular trees of infinite height (see \cite[Proposition 7.13]{PT});

\noindent $(iii)$ the volume growth of such trees is exponential (see \cite[Proposition 7.3]{PT}).

\smallskip

\noindent Different structure hypotheses on $\mathcal{T}$ give rise to polynomial volume growth (see \cite[Lemma 5.1]{FK}), and expectedly to the infimum of the spectrum of $-\D$ being equal to zero. We conjecture the same situation to apply in such case as in the Euclidean case.

\smallskip

Regarding the speed of propagation, the present results are less exhaustive than those known in the Euclidean case and in the hyperbolic case. Specifically, if $f$ is of KPP type we show that (see Theorems \ref{t1}-\ref{t4}): 

\noindent $(i)$ if $u_0$ is suitably small, for any $y\in\mathcal{T}$ and for any $c>  \hat c \,:=\, M\r_1\,\frac{b_1}{b_1-1}$
 there holds
\begin{equation*}
\lim_{t\to \infty} \sup_{\r(x)\ge  \r(y) + ct}u(x,t)\,=\, 0\,;
\end{equation*}

\noindent $(ii)$ if $u_0\not \equiv 0$, for any $y\in\mathcal{T}$ and $0<c<\check{c}\,:=\,2\sqrt{f'(0)-E_0}$
\begin{equation*}
\lim_{t\to \infty} \inf_{\r(x)<\r(y) +c t}u(x,t)\,=\, 1\,.
\end{equation*}
Expectedly, there holds $\check{c}<\hat{c}$ (see Proposition \ref{compa}).

\smallskip

In the light of the above results, also in the present case arguably there exists a speed of propagation $c_0$, and there holds $\check{c}<c_0<\hat{c}$. From a qualitative point of view, it is interesting to note that both $\check{c}$ and $\hat{c}$ increase with the increase in edge length and decrease with the increase in the number of edges emanating from any vertex
(see Remark \ref{depec}). Since $\check{c}<c_0<\hat{c}$, the same holds for $c_ 0$. Therefore, the speed of propagation is lower in regular trees with shorter edges and/or higher branching numbers. In layman's terms, the speed of propagation decreases when the number of possible path choices and/or the frequency of choice increase.

\smallskip

It is natural to wonder why the results on the speed of propagation are weaker in the present case. The above outline of results for the Euclidean case (see Subsection \ref{oueu}) pointed out the key role of travelling wave solutions, which were used to build super- and subsolutions needed to prove \eqref{propa1}-\eqref{propa2}. Essentially the same technique also worked in the hyperbolic case, the analogous solutions being the  so-called {\em horospheric waves} (see \cite{MPT}). By analogy, it seems natural to seek solutions of equation \eqref{diffeq} of the form
\begin{equation}\label{utw}
u(x,t):=p(x-ct)  
\end{equation}
for any $(x,t)\in\mathcal{T}\times\overline{\R}_+$ such that $x-ct\in\mathcal{T}$, for some {\em speed} $c\ge0$ and {\em profile} $p:\R\mapsto[0,1]$\,. 

However, travelling wave solutions on a metric graph only exist if the speed satisfies suitable restrictions dictated by the continuity of $u(\cdot,t)$ and the Kirchhoff conditions at the vertices (\cite{vb1}; see also \cite{MR}). 
It is easily seen that this is not the case for regular trees. In the parlance of \cite{vb1}, a travelling wave solution as in \eqref{utw} is called {\em isotachic}, since the speed $c$ is assumed to be the same for all edges of $\mathcal{T}$.
More generally, travelling wave solutions on a metric graph $\mathcal{G}$ are allowed to have a different speed $c=c_e$ on each edge $e\in E$. Continuity of $u(\cdot,t)$ at any interior vertex $v$ implies the equality
$$
p_{e_i}(v-c_it) \,=\, p_{e_j}(v-c_jt) \quad\text{for all $e_i,e_j \subseteq\Sig_v$\,, $i\neq j$}\,,
$$
$\Sig_v$ denoting the star centered at $v$ (see Definition \ref{deg}). Assuming differentiability of any $p_{e_i}$, from the above equality we get
$$
u_{{e_i}x}(v,t)\,=\,-c_i\,p_{e_i}'(v-c_it) \,=\,-c_j\, p_{e_j}'(v-c_jt)\,=\,u_{{e_j}x}(v,t) \quad\text{for all $i, j$ as above}\,,
$$
whence by the Kirchhoff condition
\begin{equation}\label{tolex11}
\sum_{k=1}^{d_v^+} c_k\,=\,\sum_{l=1}^{d_v^-} c_l \quad\text{for all $v\in  V\setminus \pa\mathcal{G}$}
\end{equation}
(see \eqref{lave1}). In particular, for a regular tree equality \eqref{tolex11} at any vertex $v$ with ${\rm gen}(v)=n$  reads $c=b_n c$, thus $b_n=1$, which is impossible since $b_n\ge2$ $(n\in\N)$ .

In spite of these limitations, it is possible to construct {\em symmetric} sub- and supersolutions of problem \eqref{CN} suggested by  \eqref{utw}, meaning sub- and supersolutions that depend only on the distance from the root $O$ (see Subsection \ref{symm}). The proof of Theorems \ref{t1}-\ref{t4} relies both on this construction and on suitable comparison results (see Theorems \ref{geco}-\ref{gecoxi}).

Our main results are stated in Section \ref{prore} and proved in Sections \ref{expropro}-\ref{aspropro}. Auxiliary spectral results are proven in Section \ref{spetre}, the mathematical framework and monotonicity results are presented in Section \ref{spatre}. Relevant definitions and results concerning metric graphs are collected in the Appendix.


\section{Main results}\label{prore}
\setcounter{equation}{0}

\subsection{Extinction versus propagation}

Set
$$
E_0:=\min\s(-\D)\,.
$$
Concerning propagation of fronts we can prove the following result, which is analogous to \cite[Theorem 3.3]{MPT}.

\begin{theorem}\label{t3}
Let  $(H_0)$-$(H_2)$ and $(H_4)$ be satisfied. Suppose that
\begin{equation}\label{h5}
f'(0) \,>\, E_0\,,
\end{equation}
and let $u_0\not \equiv 0$ in $\mathcal{T}$. Then for the corresponding solution of problem \eqref{CN} there holds
\begin{equation}\label{e99}
\lim_{t\to \infty} u(x,t)\,=\, 1\quad\textrm{uniformly on compact subsets of $\mathcal{T}$}.
\end{equation}
\end{theorem}

For homogeneous trees we shall prove the following companion result of Theorem \ref{t3}, concerning {\em extinction} of solutions of problem \eqref{CN} for $f$ of KPP type.
\begin{theorem}\label{t30}
Let $\mathcal{T}$ be a homogeneous tree, and let $(H_1)$, $(H_2)$ and $(H_4)$ hold. Suppose that
\begin{equation}\label{h6}
f'(0)< E_0\,,
\end{equation}
and let $u_0\not \equiv 0$ be suitably small. Then for the corresponding solution of problem \eqref{CN} there holds
\begin{equation}\label{e100}
\lim_{t\to \infty} \, u(x,t)\,=\, 0\quad\textrm{uniformly in $\mathcal{T}$}.
\end{equation}
\end{theorem}

Unsurprisingly, the same conclusion of Theorem \ref{t30} holds true if assumption $(H_4)$ is replaced by $(H_5)$.
\begin{theorem}\label{t2}
Let $\mathcal{T}$ be a homogeneous tree.
Let $(H_1)$, $(H_2)$ and $(H_5)$ be satisfied, and let $u_0\not \equiv 0$ be suitably small. Then for the corresponding solution of problem \eqref{CN} equality \eqref{e100} holds.
\end{theorem}

By  ``$u_0$ suitably small'' in Theorem \ref{t30} we mean $u_0(x)\le g(\rho(x))$ for any $x \in \mathcal{T}$, the function $g=g(\rho)$ being given by Lemma \ref{ln1}. Similarly, in Theorem \ref{t2}  $u_0$ is ``suitably small'' if $u_0(x)\le \tilde{h}(\r(x))$ $(x \in \mathcal{T})$, the function $\tilde{h}$ being defined in \eqref{detih}.


\subsection{Asymptotic speed of propagation}

If propagation occurs, an estimate from above of the speed of propagation is given by the following theorem.
\begin{theorem}\label{t1}
Let  $(H_0)$-$(H_3)$ be satisfied.  Let $u_0$ be suitably small, and let $u$ be the corresponding solution of problem \eqref{CN}. Then for any $y\in\mathcal{T}$ and for any $c> \hat c$,
\begin{equation}\label{decabis}
 \hat c \,:=\, M\r_1\,\frac{b_1}{b_1-1}\,,
\end{equation}
 there holds
\begin{equation}\label{e89}
\lim_{t\to \infty} \sup_{\r(x)\ge  \r(y) + ct}u(x,t)\,=\, 0\,.
\end{equation}
\end{theorem}
By  ``$u_0$ suitably small'' in Theorem \ref{t1} we mean $u_0(x)\le m(\rho(x))$ for any $x \in \mathcal{T}$, the function $m=m(\rho)$ being given by Lemma \ref{lesta0}.

Theorem \ref{t1} is conceptually analogous to \cite[Theorem 5.1]{AW} (observe that Theorem \ref{t1} in particular holds if $f$ is of KPP type, and equality \eqref{e89} is the counterpart of \eqref{propa1}). 
If $\mathcal{T}$ is homogeneous, there holds
\begin{equation}\label{deca}
 \hat c\,=\, M\,\frac{rb}{b-1}\,.
\end{equation}

\smallskip

The following result provides a companion estimate from below:
 \begin{theorem}\label{t4}
Let  $(H_0)$-$(H_2)$ and $(H_4)$ be satisfied, and let $\r_1<\frac{\pi}{2\sqrt{f'(0)}}\,.$ Let \eqref{h5} hold, and set
\begin{equation}\label{altrecdz}
\check{c}\,:=\,2\sqrt{f'(0)-E_0}\,.
\end{equation}
 Let $u_0\not \equiv 0$ in $\mathcal{T}$, and let $u$ be the corresponding solution of problem \eqref{CN}. Then for any $y\in\mathcal{T}$ and $0<c<\check{c}$
\begin{equation}\label{e819}
\lim_{t\to \infty} \inf_{\r(x)<\r(y) +c t}u(x,t)\,=\, 1\,.
\end{equation}
\end{theorem}

Theorems \ref{t1}-\ref{t4} provide estimates from above, respectively from below of the speed of propagation. Their
compatibility is ensured by the following result.
\begin{prop}\label{compa}
Let $\mathcal{T}$ be a homogeneous tree, and let $(H_1)$, $(H_2)$ and $(H_4)$ hold. Let \eqref{h5} be satisfied. Then there holds $\check{c}<\hat c$.
\end{prop}

\begin{remark}\label{depec}
By its very definition \eqref{decabis}, $\hat{c}$ is an increasing function of $\r_1$ and a decreasing function of $b_1$. Expectedly, $\check{c}$ has the same behaviour, since by \eqref{altrecdz} and estimate \eqref{stibe} there holds
$$
\check{c}\,\le\, 2\,\sqrt{f'(0)- \frac{(b_1-1)^2}{4(b_1\r_1)^2}}\,,
$$
and the right-hand side of the above inequality decreases both when $\r_1$ decreases and when $b_1$ increases. The expectation can be proven to hold for homogeneous trees, where
\begin{equation*}
E_0\,=\,\frac {\theta^2}{r^2}  \quad\text{with}\;\; \theta\,:=\, \arccos{\frac 1R}\,, \;\;\; R\,:=\,\frac{b+1}{2\sqrt{b}}\,>\,1\,.
\end{equation*}
(see Theorem \ref{spehone}). In fact, by \eqref{altrecdz} and the above expression of $E_0$ there holds
$$
\frac{\pa\check{c}}{\pa r}\,=\,\frac{4\theta^2}{\check{c}\,r^3}\,>\,0\,, \quad \frac{\pa\check{c}}{\pa b}\,=\,-\,\frac{\theta^2}{\check{c}\,r^2}\,\frac{1}{R\sqrt{R^2-1}}\,\frac{b-1}{b}\,<\,0\,.
$$
\end{remark}


\section{Spectrum of the Neumann Laplacian on regular trees}\label{spetre}

Set
$$
L(\mathcal{T})\,:=\int_0^{\infty}\!\!\frac{d\r}{\b(\r)}\,, \qquad B(\mathcal{T})\; := \; \sup_{t>0}\,\left(\int_0^t \b(\r)\, d\r\int_t^\infty \frac{d\r}{\b(\r)}\right)\,,
$$
where $\b:\overline{\R}_+\mapsto \N$ is the branching function of $\mathcal{T}$ (see \eqref{bratree}).
Then the following holds:
\begin{theorem}\label{tre00}
Let  $(H_0)$ be satisfied. Then $E_0>0$.
\end{theorem}
To prove Theorem \ref{tre00} the following lemma is needed.

\begin{lem}\label{LBfinite}
Let  $(H_0)$ be satisfied. Then $L(\mathcal{T})<\infty$, $B(\mathcal{T})<\infty$\,.
\end{lem}
\begin{proof} Set $A_m:= b_0b_1\dots b_m$ $(m\in\N\cup\{0\})$, thus by $(H_0)$-$(ii)$ there holds $A_m\ge b_1^m$. Set $I_k:=(\r_{k-1},\r_k)$ $(k\in\N)$. Then from \eqref{bratree} we get
\begin{equation*}
L(\mathcal{T})\,=\sum_{k=0}^\infty\,\int_{I_k}\frac{dr}{\b(\r)}\,=\,
\sum_{k=1}^\infty \frac{\r_k-\r_{k-1}}{A_{k-1}}\,\le\, \r_1\sum_{k=1}^\infty \frac{1}{b_1^{k-1}} \,=\,\frac{b_1\r_1}{b_1-1}\,.
\end{equation*}

Now observe that $A_n=b_0\dots b_m b_{m+1}\dots b_n=A_m b_{m+1}\dots b_n$ if $m<n$\,. Let $n\in\N$ be fixed and $t\in I_{n+1}\equiv(\r_n,\r_{n+1})$\,. Then we get plainly
\begin{eqnarray*}
&&\int_0^t \b(\r)\, d\r\int_t^\infty \frac{d\r}{\b(\r)}\;\le\,\int_0^{\r_{n+1}}\!\! \b(\r)\, d\r\int_{\r_n}^\infty \frac{d\r}{\b(\r)}\,=\\
&=&\left(\sum_{k=1}^{n+1} A_{k-1}(\r_k-\r_{k-1})\right)\left(\sum_{k=n+1}^\infty \frac{\r_k-\r_{k-1}}{A_{k-1}}\right)\;\le\;
 \r_1^2\, \left(\sum_{l=0}^n \frac{A_l}{A_n}\right)\left(\sum_{l=0}^\infty \frac{1}{A_l}\right)\,\le\\
 &\le&\r_1^2 \left(\sum_{m=0}^n \frac{1}{b_1^m}\right)\left(\sum_{m=0}^\infty \frac{1}{b_1^m}\right)\;\le\; \r_1^2 \,\frac{b_1}{(b_1-1)^2\,}\left(b_1- \frac{1}{b_1^n}\right)\qquad (t\in I_{n+1}) \,.
\end{eqnarray*}
Then there holds
\begin{equation}\label{sdue}
B(\mathcal{T})\,\le\,  \left(\frac{b_1\r_1}{b_1-1}\right)^2\,,
\end{equation}
whence the result follows.
\end{proof}
Observe that, if $\mathcal{T}$ is a homogeneous tree with $b_n=b$, $\r_n=nr$ $(n\in\N, r>0)$, the proof of Lemma \ref{LBfinite} shows that $L(\mathcal{T})=\frac{br}{b-1}$.

\smallskip

\noindent {\em Proof of Theorem \ref{tre00}.}
It is known that, if $\mathcal{T}$ is a regular tree of infinite height,  the following statements are equivalent (see \cite[Theorem 5.2 and Subsection 7.1]{So2}):

\noindent  $(i)$ $E_0>0\,;$

\noindent $(ii)$ there holds $L(\mathcal{T})\,<\, \infty$\,, $B(\mathcal{T})\,<\,\infty$\,.

\noindent Then by Lemma \ref{LBfinite} the conclusion follows.
\hfill$\square$
\begin{remark}
An additional result of \cite[Theorem 5.2]{So2} is that, if $E_0>0$, there holds
\begin{equation*}
 \frac{1}{4B(\mathcal{T})} \,\le\,   E_0\, \le\, \frac{1}{B(\mathcal{T})} \,.
\end{equation*}
From the above inequality and \eqref{sdue}, under assumption $(H_0)$-$(ii)$ we get the estimate
\begin{equation}\label{stibe}
E_0\,\ge\,  \frac{(b_1-1)^2}{4(b_1\r_1)^2} \,.
\end{equation}
\end{remark}

Remarkably, if $\mathcal{T}$ is homogeneous, $E_0$ can be explicitly calculated:

\begin{theorem}\label{spehone}
Let  $\mathcal{T}$ be homogeneous with  $b_n=b$ and edge length $r>0$. Then there holds
\begin{equation}\label{isp}
E_0\,=\,\frac {\theta^2}{r^2} \,, \quad\text{where}\;\; \theta\,:=\, \arccos{\frac 1R}\,, \;\;\; R\,:=\,\frac{b+1}{2\sqrt{b}}\,>\,1\,.
\end{equation}
\end{theorem}
\begin{proof} By  \cite[Corollary 3.7]{So2} there holds $E_0\,=\, \min \s(A)$, $A$ being the operator associated to the form $\mathfrak{a}_0$ with domain in $L^2(\R_+)$ defined by \eqref{mnbbis}. Hence the conclusion follows from the following

 \smallskip

 \noindent {\em Claim:} There holds
\begin{equation}\label{spettro}
\s(A)
\,=\, \bigcup_{l\in\N} \left[\left(\frac{\pi(l-1)+\theta}{r}\right)^2\!\!,\,\left(\frac{\pi l-\theta}{r}\right)^2\right] \,.
\end{equation}

The same argument used in the proof of \cite[Theorem 3.2]{SS} shows that the essential spectrum of the operator $A$ coincides with the right-hand side of \eqref{spettro}. Then the Claim will follow if we prove that the discrete spectrum of $A$ is empty, namely that no isolated eigenvalues of $A$ of finite multiplicity exist  in $\overline{\R}_+\setminus \left(\bigcup_{l\in\N} \left[\left(\frac{\pi(l-1)+\theta}{r}\right)^2\!\!,\,\left(\frac{\pi l-\theta}{r}\right)^2\right]\right)$\,.

To this purpose, observe preliminarily by direct inspection that $\l=0$ is not an eigenvalue. Further, let us seek  $\z\in L^2(\R_+)$ satisfying
$$
-A\z +\l \z=0 \quad\text{a.e. in $\R_+$} \qquad (\l>0)\,.
$$
Then $\z$ must belong to the domain of $A$, thus by Proposition \ref{frado} is of the form $\z=
\sum_{n=1}^\infty \z_n \chi_{_{I_n}}$  and satisfies the following requirements:
\begin{subequations}\label{auto}
\begin{equation}\label{auto2}
 \z_n''\,+\,\l\z_n\,=\,0 \quad\textrm{ in $I_n$ for all $n\in\N$\,,}
\end{equation}
\begin{equation}\label{auto3}
\z_{n+1}(\r_n)\,=\, \sqrt{b}\,\z_n(\r_n)\,,\quad \z_n'(\r_n^-)\,=\, \sqrt{b}\,\z_{n+1}'(\r_n^+)\quad\text{for all $n\in\N$}\,,
\end{equation}
\begin{equation}\label{auto4}
\z_1'(0^+)\,=\,0\,.
\end{equation}
\end{subequations}

Consider the quadratic equation
\begin{equation}\label{fqe}
\a^2 -2R\cos{(\sqrt{\l}r)}\a+ 1=0\,,
\end{equation}
where $R\,:=\,\frac{b+1}{2\sqrt{b}}$. If $R|\cos{(\sqrt{\l}r)}|>1$, equation \eqref{fqe} has two real roots
\begin{equation}\label{roqe}
\a_\pm \,=\, R\cos{(\sqrt{\l}r)} \,\pm\,\sqrt{R^2\cos^2{(\sqrt{\l}r)} -1}
\end{equation}
(observe that $R\,|\cos{(\sqrt{\l}r)}|\le1$ exactly when $\l$ belongs to the right-hand side of \eqref{spettro}).
Following \cite{SS}, it is easily checked that for any $n\in\N$ the functions
\begin{equation}\label{exqn}
\z_{n,\pm}(\r)\,:=\, \a_\pm^{n-1}\left\{\sqrt{b}\sin{[\sqrt{\l}(nr-\r)]} \,+\, \a_\pm\sin{[\sqrt{\l}(\r-(n-1)r)]}\right\} \qquad(\r\in I_n=((n-1)r,nr))
\end{equation}
satisfy \eqref{auto2}-\eqref{auto3}, and their Wronskian is
$$
\z_{n,+}(\r)\z_{n,-}'(\r) \,- \, \z_{n,+}'(\r)\z_{n,-}(\r)\,=\,\sqrt{b} (\a_- -\a_+)\sqrt{\l}\sin{(\sqrt{\l} r)} \,.
$$

Therefore, if $\l\neq\left(\frac{m\pi}{r}\right)^2$ $(m\in\N)$, the functions defined in \eqref{exqn} are linearly independent, and every $\z$ such that \eqref{auto2}-\eqref{auto3} hold is of the form
$$
\z=c_+\z_+ +  c_-\z_- \quad\text{with $c_\pm\in\R$, \; $\z_\pm:=
\sum_{n=1}^\infty \z_{n,\pm} \chi_{_{I_n}}$\,.}
$$
Observe that by \eqref{exqn} $\z$ does not belong to $L^2(\R_+)$, since $\a_+>1$ when $ R\cos{(\sqrt{\l}r)}>1$ and $\a_-<-1$ when $ R\cos{(\sqrt{\l}r)}<-1$.

Imposing the Neumann condition at $\r=0$ plainly gives
$$
c_-=-\frac{\sqrt{b}\cos{(\sqrt{\l}r)} -\a_+}{\sqrt{b}\cos{(\sqrt{\l}r)} -\a_-}\,c_+\,.
$$
Hence every function $\z$ satisfying \eqref{auto} is proportional to the function
$$
\frac{\sqrt{b}\cos{(\sqrt{\l}r)}\left[\z_+(\r) -\z_-(\r)\right] -\left[q_-\z_+(\r) -q_+\z_-(\r)\right]}{\sqrt{b}\cos{(\sqrt{\l}r)} -\a_-}\,,
$$ that does not belong to $L^2(\R_+)$\,.

To sum up,  we proved that no eigenfunctions of $A$ exist if $R\,|\cos{(\sqrt{\l}r)}|>1$ and $\l\neq\left(\frac{m\pi}{r}\right)^2$
$(m\in\N)$\,. If for some $m\in\N$ there holds $\l=\left(\frac{m\pi}{r}\right)^2$,  it is easily seen that every $\z$ satisfying \eqref{auto} is proportional to the function
$$
\sum_{n=1}^\infty b^{n/2}\cos{\left(\frac{m\pi\r}{r} \right)} \chi_{_{I_n}}(\r)\,,
$$
which does not belong to $L^2(\R_+)$ either. This completes the proof.
\end{proof}

We finish this section by proving Proposition \ref{compa}.

\smallskip

\noindent {\em Proof of Proposition \ref{compa}.}
In view of \eqref{isp}, we have that
\begin{equation}\label{risp}
\check{c}\,<\,\hat c \quad \Leftrightarrow \quad 2\sqrt{Mr^2-\theta^2}\,<\, Mr^2\,\frac{b}{b-1}\,,
\end{equation}
with $\theta\equiv\theta(b):= \arccos{\frac{2\sqrt{b}}{b+1}}$\,. It is easily seen that the function $\theta(b)$ is increasing and there holds $\theta(2)<\theta(b)<\theta(0)=\frac\pi 2$, with $\theta(2)=\arccos{\frac{2\sqrt{2}}{3}}\sim \frac{\pi}{12}<1$\,.

 Since $\frac{b}{b-1}\ge1$ for any $b\in[2,\infty)$, the second inequality in \eqref{risp} is satisfied if
$$
2\sqrt{Mr^2-\theta^2}\,<\, Mr^2 \quad \Leftrightarrow \quad (Mr^2 )^2-4Mr^2+4\theta^2\,>\, 0\,.
$$
An elementary analysis shows that the latter inequality is satisfied for any value of $Mr^2>0$ if $\theta>1$, or for $Mr^2>
2+2\sqrt{1-\theta^2}$ if $\theta\le1$. Since  $Mr^2>\theta^2$ by \eqref{h5} and \eqref{isp}, and $\theta^2>2+2\sqrt{1-\theta^2}$, the second inequality in \eqref{risp} is satisfied in this case, too. Hence the result follows.
\hfill$\square$


\section{Semilinear parabolic equations on regular trees}\label{spatre}
\setcounter{equation}{0}

\subsection{Concept of solution and well-posedness}\label{spe}

Solutions of problem \eqref{CN} are meant in the following sense.
\begin{definition}\label{deso}
Let $(H_1)$ hold. A  function $u\in C(\overline{\R}_+;L^2(\mathcal{T}))\cap C^1(\R_+;L^2(\mathcal{T}))$, $0\le  u(\cdot,t) \le  1$ a.e. in $\mathcal{T}$ for any $t\in \overline{\R}_+$, such that $\D u \in C(\R_+;L^2(\mathcal{T}))$, is called a {\em solution} of problem \eqref{CN}  if
$$
\left\{
\begin{array}{ll}
u_t(\cdot,t) = (\Delta u)(\cdot,t)+ f(u(\cdot,t))  \quad\text{a.e. in $\mathcal{T}$ for any $t\in \R_+$\,,}
\\& \\
u(\cdot,0)= u_0\quad\text{a.e. in $\mathcal{T}$}.
\end{array}
\right.
 $$
\end{definition}

\begin{remark}\label{condege}
$(i)$ By Definition \ref{deso}, for any $t\in\R_+$ $u(\cdot,t)$ belongs to the domain of the Neumann Laplacian. Then by \eqref{lave} for all $t\in \R_+$ there holds $u(\cdot,t)=\bigoplus_{e\in E} u_e(\cdot,t)$, and
\begin{subequations}\label{solge}
\begin{equation}\label{solge0}
 u(\cdot,t)\in H^1(\mathcal{T})\,,\quad  u_e(\cdot,t)\in H^2(I_e) \,\textrm{ for all $e\in E$}, \;
 \quad  \sum_{e\in E}\int_0^{l_e}
 | u_{exx}(\cdot,t)|^2\,dx<\infty \,,
\end{equation}
\begin{equation}\label{solge1}
u_{et}(\cdot,t)  \,=\,  u_{exx}(\cdot,t) + f(u_e(\cdot,t))  \quad\text{a.e. in $I_e$ for all $e\in E$\,,}
\end{equation}
\begin{equation}\label{solge2}
u_x (O^+,t)  \,=\,0\,, \qquad \sum_{e\ni v}\frac{d u_e}{d\n}(v,t)=0 \quad\text{for all $v\in \mathcal{T}\setminus \{O\}$}\,.
\end{equation}
\end{subequations}
By \eqref{solge0} and embedding results, for any $t\in \R_+$:
$(a)$ $u(\cdot,t)$ is continuous in $\mathcal{T}$, thus at any vertex $v\in \mathcal{T}\setminus \{O\}$; $(b)$  there holds $u_e(\cdot,t)\in C^1(\overline{I}_e)$ $(e\in E)$.

\noindent$ (ii)$ For any vertex $v$ of the $n$-th generation, thus belonging to $\mathcal{T}\setminus O$, the inbound star $\Sig_v^+$  centered at $v$ consists of one edge $e$ (i.e., $d_v^+=1$), whereas the outbound star $\Sig_v^-$  consists of $b_n$ edges $e_l$ (i.e., $d_v^-=b_n$; see Definition \ref{deg}).  Therefore, recalling that the maps $ i,j:E\mapsto V$  define the initial point, respectively the  final point of an edge (see Definition \ref{meg}), the second equality in \eqref{solge2} reads \begin{equation}\label{tolex111}
u_{ex}(j(e),t)\,=\,\sum_{l=1}^{{b_n}} u_{{e_l}x}(i(e_l),t) \qquad (t\in\R_+,\,n\in\N) \,.
\end{equation}
\end{remark}

The following result, which ensures well-posedness  of problem \eqref{CN}, follows by standard results of semigroup theory (e.g., see \cite[Theorem 4.4]{Y}).
\begin{theorem}\label{wpp}
Let  $(H_1)$-$(H_2)$ hold. Then there exists a unique solution of problem \eqref{CN}.
\end{theorem}

\begin{definition}\label{destage}
A function $q:\mathcal{T}\mapsto [0,1]$, $q\in D(\D)$ is a {\em stationary solution} of  equation \eqref{diffeq} if
\begin{equation}\label{stage}
\D q + f(q) =0  \quad\text{a.e. in $\mathcal{T}$.}
\end{equation}
\end{definition}
\noindent Equivalently, a stationary solution $q=\bigoplus_{e\in E} q_e$ of \eqref{diffeq} satisfies the following:
\begin{subequations}\label{proge}
\begin{equation}\label{proge1}
q\in H^1(\mathcal{T})\,,\quad  q_e\in H^2(I_e) \,\textrm{ for all $e\in E$}, \;
 \quad  \sum_{e\in E}\int_0^{l_e}| q_e''|^2\,dx<\infty \,,
\end{equation}
\begin{equation}\label{proge2}
 q_e''+f(q_e)\,=\,0  \quad\text{a.e. in $I_e$ for all $e\in E$\,,}
\end{equation}
\begin{equation}\label{proge3}
q'(O^+)\, =\,0\,, \qquad \sum_{e\ni v}\frac{d q_e}{d\n}(v)=0 \quad\text{for all $v\in \mathcal{T}\setminus \{O\}$}\,.
\end{equation}
\end{subequations}

\begin{remark}\label{uniso=1}
The function $q\equiv1$ is not a stationary solution of \eqref{diffeq} in the sense of Definition \ref{destage}, since it does not belong to $L^2(\mathcal{T})$. However, it is a {\em weak stationary solution} of \eqref{diffeq}, in the sense that
\begin{equation}\label{stawe}
\int_\mathcal{T} \{q\D\eta + f(q)\eta\}\,d\m \,=\,0\quad\text{ for all $\eta\in C^2_c(\mathcal{T})$}\,.
\end{equation}
Moreover, it is the only nontrivial solution of  \eqref{stawe} (in this connection, see \cite[Theorem 3]{vb1}). In fact, if $q=\bigoplus_{e\in E} q_e$ satisfies \eqref{stawe}, by the arbitrariness of $\eta$ and regularity results $q_e$ is a classical solution of problem \eqref{proge2} for all $e\in E$, and \eqref{proge3} holds. In view of assumption $(H_2)$, the maximum principle and \eqref{proge3}, it is easily seen that no such a solution $q$ of \eqref{CN}, $q:\mathcal{T}\mapsto(0,1]$ nonconstant, exists. Hence the claim follows.
\end{remark}


\subsection{Monotonicity results}\label{wepomo}

Super- and subsolutions of problem \eqref{CN} are defined as follows.
\begin{definition}\label{evoss}
A function $\overline{u}=\bigoplus_{e\in E} \overline{u}_e\in C(\overline{\R}_+;L^2(\mathcal{T}))\cap C^1(\R_+;L^2(\mathcal{T})))\cap C(\R_+;H^1(\mathcal{T}))$, $0\le  \overline{u}(\cdot,t) \le  1$ a.e. in $\mathcal{T}$ for any $t\in \overline{\R}_+$, is a {\em supersolution} of problem \eqref{CN} if:
\begin{subequations}\label{sottop}
\begin{equation}\label{sottop0}
 \overline{u}_e(\cdot,t)\in H^2(I_e) \,\textrm{ for all $e\in E$}, \;
 \quad  \sum_{e\in E}\int_0^{l_e}| \overline{u}_{exx}(\cdot,t)|^2\,dx<\infty\quad\text{for any $t\in \R_+$\,,}
\end{equation}
\begin{equation}\label{sottop1}
\overline{u}_{et}(\cdot,t)\,\ge\,  \overline{u}_{exx}(\cdot,t)+  f(\overline{u}_e(\cdot,t))  \quad\text{a.e. in $I_e$ for any $e\in E$ and $t\in \R_+$\,,}
\end{equation}
\begin{equation}\label{solprop22}
\overline{u}(\cdot,0)\,\ge\, u_0 \;\;\text{a.e. in $\mathcal{T}$}\,,
\end{equation}
\begin{equation}\label{solprop23}
\overline{u}_x(O^+,t)\le0\,, \qquad \sum_{e\ni v}\frac{d \overline{u}_e}{d\n}(v,t)\,\ge\,0 \quad\text{for all $v\in \mathcal{T}\setminus \{O\}$ and $t\in\R_+$}\,.
\end{equation}
\end{subequations}
A {\em subsolution} $\underline{u}$ of \eqref{CN} is defined by reversing inequalities in \eqref{sottop1}-\eqref{solprop23}.
\end{definition}

\begin{remark}\label{limiss}
$(i)$ Since $H^1(\mathcal{T})\subseteq C(\mathcal{T})$ and $H^2(I_e)\subseteq C^1(\overline{I}_e)$, by \eqref{sottop0} for all $t\in\R_+$:

\noindent $(a)$ both $\underline{u}(\cdot,t)$ and $\overline{u}(\cdot,t)$ are continuous in $\mathcal{T}$, thus at any vertex $v\in \setminus \{O\}$;

\noindent $(b)$ there holds $\underline{u}_e(\cdot,t), \overline{u}_e(\cdot,t)\in C^1(\overline{I}_e)$, thus in particular $\underline{u}_x(\cdot,t)$ and $\overline{u}_x(\cdot,t)$ are continuous (from the right) at the root $O$.

\noindent $(ii)$ In analogy with \eqref{tolex111}, the second inequality in \eqref{solprop23} reads
\begin{equation}\label{stolex}
u_{ex}(j(e),t)\,\ge\,\sum_{l=1}^{{b_n}} u_{{e_l}x}(i(e_l),t) \qquad (t\in\R_+,\,n\in\N) \,.
\end{equation}
\end{remark}

\begin{definition}\label{soso}
A function $\overline{q}:\mathcal{T}\mapsto [0,1]$, $\overline{q}=\bigoplus_{e\in E} \overline{q}_e\in H^1(\mathcal{T})$ is a {\em stationary supersolution} of equation \eqref{diffeq} if
\begin{subequations}\label{sproge}
\begin{equation}\label{sproge1}
\overline{q}_e\in H^2(I_e) \,\textrm{ for all $e\in E$}, \;
 \quad  \sum_{e\in E}\int_0^{l_e}| \overline{q}_e''(x)|^2\,dx<\infty \,,
\end{equation}
\begin{equation}\label{sproge2}
 \overline{q}_e''+f(\overline{q}_e)\,\le\,0  \quad\text{a.e. in $I_e$ for all $e\in E$\,,}
\end{equation}
\begin{equation}\label{sproge3}
 \overline{q}'(O^+)\le0\,,\qquad \sum_{e\ni v}\frac{d  \overline{q}_e}{d\n}(v)\,\ge\,0 \quad\text{for all $v\in \mathcal{T}\setminus \{O\}$}\,.
\end{equation}
\end{subequations}
A {\em stationary subsolution} $\underline{q}$ is defined by reversing inequalities in \eqref{sproge2}-\eqref{sproge3}.
\end{definition}

A different concept of super- and subsolutions is the content of the following definition.
\begin{definition}\label{sosol2}
A function $ \overline{u}\in C(\mathcal{T}\times \overline{\R}_+)\cap L^\infty(\mathcal{T}\times \R_+), \overline{u}(\cdot,t)=\bigoplus_{e\in E} \overline{u}_e(\cdot,t)$ is  a {\em bounded supersolution}  of problem \eqref{CN} if:

\noindent $(i)$ for all $e\in E$ there holds
\begin{subequations}\label{mnbv}
\begin{equation}\label{mnbv0}
\text{ $\overline{u}_e(\cdot, t)\in C^2(I_e)\cap C^1(\bar{I}_e)$ for any $ t\in \R_+\,,$}
\end{equation}
\begin{equation}\label{mnbv1}
\text{ $\overline{u}_e(x, \cdot)\in C^1(\R_+)$ for any $x\in \bar{I}_e \,;$}
\end{equation}
\end{subequations}

\noindent $(ii)$ inequalities \eqref{sottop1}-\eqref{solprop23} hold.

\noindent A {\em bounded subsolution} of \eqref{CN} is similarly defined, with reverse inequalities in \eqref{sottop1}-\eqref{solprop23}.
\end{definition}

We refer the reader to \cite{PTmo} for the proof of the following results.
\begin{theorem}\label{geco}
Let $(H_1)$-$(H_2)$ hold, and let $\underline{u}, \,\overline{u}$ be a sub-, respectively a supersolution of problem \eqref{CN}. Then there holds $\underline{u}(x,t) \le \overline{u}(x,t)$ for any $(x,t)\in \mathcal{T}\times \overline{\R}_+$.
\end{theorem}

\begin{theorem}\label{decre}
Let $(H_1)$-$(H_2)$ hold. Let $\underline{q},\,\overline{q}$ be a stationary sub-, respectively supersolution of  \eqref{diffeq}, such that $\underline{q}\le\overline{q}$ a.e. in $\mathcal{T}$. Let $u_1$, $u_2$ be the solutions of \eqref{CN} with $u_0=\underline{q},$ respectively $u_0=\overline{q}$. Then:

\noindent $(i)$ $u_1$, $u_2$ are global, and for all $t\in\R_+$ there holds $\underline{q}\le u_1(\cdot,t)\le u_2(\cdot,t)\le\overline{q}$ in $\mathcal{T}$;

\noindent $(ii)$ if $\underline{q}\le u_0\le\overline{q}$, the solution $u$ of problem \eqref{CN}  is global, and for all $t\in\R_+$ there holds $u_1(\cdot,t)\le u(\cdot,t)\le u_2(\cdot,t)$ in $\mathcal{T}$;

\noindent $(iii)$ the map $t\mapsto u_1(\cdot,t)$ is nondecreasing, $t\mapsto u_2(\cdot,t)$ is nonincreasing, and the pointwise limits $\hat{u}_1:=\lim_{t\to\infty}u_1(\cdot,t)$, $\hat{u}_2:=\lim_{t\to\infty}u_2(\cdot,t)$ are stationary solutions of \eqref{diffeq}. Moreover, there holds $\hat{u}_1\le \hat{u}_2$\,, and every stationary solution of \eqref{diffeq} with $\underline{q}\le \tilde{u}\le\overline{q}$ satisfies $\hat{u}_1\le\tilde{u}\le\hat{u}_2$\,.
\end{theorem}

\begin{theorem}\label{comparison}
Let $(H_1)$-$(H_2)$ hold, and let $\underline u$, $\overline u$ be a  bounded subsolution, respectively a  bounded supersolution of \eqref{CN}. Then there holds  $\underline{u}(x,t) \le \overline{u}(x,t)$ for any $(x,t)\in \mathcal{T}\times \overline{\R}_+$.
\end{theorem}

A companion result of Theorem \ref{decre} follows from Theorem \ref{comparison}; we leave its formulation to the reader.

\begin{remark}\label{wefo} 
Arguing as in \cite{PTmo}, it is easily seen that Theorem \ref{comparison} still holds for {\em weak} bounded super- and subsolutions, in the sense of the following definition.

\begin{definition}
A function $ \overline{u}\in  C(\overline{\R}_+; H_{loc}^1(\mathcal T))\cap L^\infty(\mathcal{T}\times \R_+), \overline{u}(\cdot,t)=\bigoplus_{e\in E} \overline{u}_e(\cdot,t)$ is a {\em weak bounded supersolution}  of problem \eqref{CN} if:

\noindent $(i)$ for any $T>0$ and $\varphi\in C^1([0,T]; C_c^1(\mathcal T))$ such that $\varphi\geq 0$, $\varphi(\cdot,T)=0$ in $\mathcal T$ there holds
\begin{equation}\label{e1000}
\int_0^T\!\!\!\int_{\mathcal T} \overline{u} \varphi_t \,dx dt \;\leq\, \int_0^T\!\!\!\int_{\mathcal T}\left[\overline{u}_x\varphi_x - f(\overline{u}) \varphi \right]dx dt \,-\, \int_{\mathcal T} u_0 \varphi(x,0)\, dx \,;
\end{equation}

\noindent $(ii)$ the derivatives in \eqref{solprop23} exist, and inequalities \eqref{solprop23} are satisfied. 

\noindent A {\em weak bounded subsolution} of \eqref{CN} is similarly defined, with reverse inequalities in \eqref{solprop23} and in \eqref{e1000}.
\end{definition}
\end{remark}

\begin{remark}\label{dinico}
It is worth mentioning that the convergences $u_i(\cdot,t)\to\hat{u}_i$ $(i=1,2)$ in Theorem \ref{decre}-$(iii)$ are uniform in any compact $K\subset\mathcal{T}$. In fact, let  a compact $K\subset\mathcal{T}$ and $\t>0$ be fixed, and let $\varphi\in C_c^\infty(K\times (0,\t))$, $\varphi \ge0$\,. Since the map $t\mapsto u_2(\cdot,t)$ is nonincreasing, there holds
\begin{equation}\label{regev}
\int_0^\t\!\!\int_K\left[ u_2\, \D \varphi +f(u_2)\varphi\right]\,dx dt \,\le\,0\,.
\end{equation}
Choosing $\varphi(x,t)=\phi(x)\theta(t)$ with $\phi\in C_c^\infty(K)$, $\theta\in C_c^\infty(0,\t)$, $\phi\ge0$, $\theta\ge0$, from \eqref{regev} and assumption $(H_2)$ we obtain that $\D u_2 (\cdot, t)+ f(u_2(\cdot,t))\le0$. Then $\D u_2 (\cdot, t)$ is a negative distribution, hence a measure on $K$ for any $t\in(0,\t)$. Then by embedding results and the arbitrariness of $\t$ there holds $u_2(\cdot,t)\in W^{1,\infty}(K) \subseteq C(K)$ for any $t>0$. On the other hand, since $\hat{u}_2$ is a  stationary solution of \eqref{diffeq}, there holds $\hat{u}_2 \in C(K)$ (see  \eqref{proge1}). Then by the Dini Theorem the claim concerning $u_2$ follows; a similar argument applies to $u_1$\,.
\end{remark}


\subsection{Symmetric super- and subsolutions}\label{symm}

It is natural to seek {\em symmetric} solutions and supersolutions or subsolutions of problem \eqref{CN}  (see Definition \ref{symm}). Symmetric solutions (respectively supersolutions, subsolutions) are of the form $u(x,t)=z(\r(x),t)$ (respectively $\overline{u}(x,t)=\overline{z}(\r(x),t)$, $\underline{u}(x,t)=\underline{z}(\r(x),t)$), where $z,\overline{z},\,\underline{z}:\overline{\R}_+\times\overline{\R}_+\mapsto \R$ and $\r=\r(x):=d(x,O)$ $(x\in\mathcal{T})$.

In the light of of the above remarks, it is immediately seen that $u$ is a symmetric solution of problem \eqref{CN}  in $\overline{\R}_+$ if and only if (see \eqref{dl2b} and following remarks):
\begin{subequations}\label{sotto}
\begin{equation}\label{sotto0}
z\in C(\overline{\R}_+;L^2(\R_+;\b))\cap\, C^1(\R_+;L^2(\R_+;\b))\cap C(\R_+;H^1(\R_+;\b))\,,
\end{equation}
\begin{equation}\label{sotto10}
\text{$0\le  z(\cdot,t) \le  1$ a.e. in $\R_+$\,, \quad$z(\r,t)=\sum_{n=1}^\infty z_n(\r,t)\chi_{_{I_n}}(\r)$ \;\;for any $t\in \R_+$\,,}
\end{equation}
\begin{equation}\label{sotto20}
 z_n(\cdot,t)\in H^2(I_n)\;\textrm{ for all $n\in\N$}\,,\quad \sum_{n\in\N}\int_{I_n} \left|z_{n\r\r} (\r,t) \right|^2\b(\r)  d\r\,<\,\infty \quad\text{for any $t\in \R_+$\,,}
\end{equation}
\begin{equation}\label{sotto200}
z_t(\cdot,t)\,=\,  z_{\r\r}(\cdot,t) + f(z(\cdot,t))  \quad\text{a.e. in $\R_+$ for any $t\in \R_+$\,,}
\end{equation}
\begin{equation}\label{sotto3}
z(\r_n^-,t)\,=\,z(\r_n^+,t)\,,\quad z_{\r}(\r_n^-,t)\,=\, b_n\,z_{\r}(\r_n^+,t)\quad\text{for all $n\in\N$\,,}
\end{equation}
\begin{equation}\label{sotto2}
z_{\r}(0^+\!,t) \,=\,0 \;\text{ for any $t\in \R_+$\,,} \quad z_0(\r(x))\,:=\, z(\r(x),0)\,=\, u_0(x) \;\;\text{for a.e. $x\in\mathcal{T}$}\,.
\end{equation}
\end{subequations}
Observe that the first equality in \eqref{sotto3} follows from the contiunity of $z(\cdot,t)$ in $\R_+$ (see Remark \ref{marko}).
Symmetric super- and subsolutions of \eqref{CN} have properties \eqref{sotto0}-\eqref{sotto20} and satisfy \eqref{sotto200}-\eqref{sotto3} with proper inequalities (see Definition \ref{evoss}).

Similar considerations hold for symmetric stationary sub- and supersolutions of equation \eqref{diffeq} - namely, for sub- and supersolutions of equation \eqref{diffeq} of the form  $\underline{q}(x)= \underline{p}(\r(x))$, $\overline{q}(x)= \overline{p}(\r(x))$ $(x\in\mathcal{T})$. For instance, in view of Definition \ref{soso}, $\overline{q}$ is a symmetric stationary supersolution of equation \eqref{diffeq} if and only if $\overline{p}$ belongs to $H^1(\R_+;\b)$,
$\overline{p}=\sum_{n=1}^\infty \overline{p}_n\chi_{_{I_n}}$, $0\le\overline{p}\le1$ in $\R_+$\,,  and the following holds:
\begin{subequations}\label{sproge}
\begin{equation}\label{sproge1}
\overline{p}_n\in H^2(I_n) \,\textrm{ for all $n\in \N$}, \;
 \quad  \sum_{n\in \N}\int_{l_n}| \overline{p}_n''(\r)|^2\b(\r)\,d\r<\infty \,,
\end{equation}
\begin{equation}\label{sproge2}
 \overline{p}''+f(\overline{p})\,\le\,0  \quad\text{a.e. in $\R_+$\,,}
\end{equation}
\begin{equation}\label{sproge3}
\overline{p}'(0^+)\le0\,, \quad \overline{p}(\r_n^-)\,=\, \overline{p}(\r_n^+)\,, \quad \overline{p}'(\r_n^-)\,\ge\, b_n\,\overline{p}'(\r_n^+)\quad\text{for all $n\in\N$}\,,
\end{equation}
\end{subequations}

Similar remarks and results hold for symmetric bounded super- and subsolutions of problem \eqref{CN}; we leave their formulation to the reader. 

\smallskip

We shall be interested in studying symmetric solutions of problem \eqref{CN} restricted to the family of half-lines $H_n:=\{(\r,t)\,|\, \r=\r_n+ct, t\ge0\,; n\in\N\}$, thus to the domain $D:=\{(\r,t)\,|\, \r\ge ct, t\ge0 \}$ (see the proofs of Theorems \ref{t1}-\ref{t4}). Set $\xi:=\r-ct$ with $\r\ge ct$, $t\in\R_+$\,. By the change of variables $(\r,t) \mapsto (\xi,t)$ the domain $D$ is mapped onto the first quadrant $Q_+:=\{(\xi,t)\,|\, \xi\ge0, t\ge0\}$, and each $H_n$ onto the half-line $\{(\xi,t)\,|\, \xi=\r_n,\, t\ge0\}$. Let $J_n(t):=(\xi_{n-1}(t),\xi_n(t))$ with $\xi_n(t):=\max\{\r_n-ct,0\}$ ($n\in\N\cup\{0\},t\in\overline{\R}_+$); observe that $\xi_0\equiv0$, and $J_n(t)=\emptyset$ for any $n=1,\dots,k$ if $t\ge t_k:=\frac{\r_k}{c}$ $(k\in\N)$. Also set $t_0:=0$,  $\chi_\emptyset:=0$\,.

Plainly, $z$ satisfies \eqref{sotto0}-\eqref{sotto3} if and only if $v=v(\xi,t):=z(\xi+ct,t)$ satisfies the following:

\begin{subequations}\label{sottsass}
\begin{equation}\label{sottsass01}
v\in C(\overline{\R}_+;L^2(\R_+;\b))\cap\, C^1(\R_+;L^2(\R_+;\b))\cap C(\R_+;H^1(\R_+;\b))\,,
\end{equation}
\begin{equation}\label{sottsass02}
0\le  v(\cdot,t) \le  1\;\;\text{a.e. in $\R_+$}\,, \quad v(\xi,t)=\sum_{n=1}^\infty v_n(\xi,t)\chi_{_{J_n(t)}}(\xi) \qquad (\xi\in J_n(t), t\ge0)\,,
\end{equation}
\begin{equation}\label{sottsass20}
 v_n(\cdot,t)\in H^2(J_n(t))\;\textrm{ for all $n\in\N$}\,,\quad \sum_{n=1}^\infty\int_{J_n(t)} \left|v_{n\xi\xi} (\xi,t) \right|^2\b(\xi)  d\xi\,<\,\infty \quad\text{for any $t\in \R_+$\,,}
\end{equation}
\begin{equation}\label{sottsass1}
v_t(\cdot,t)\,=\, v_{\xi\xi}(\cdot,t) + cv_\xi(\cdot,t) + f(v(\cdot,t))  \quad\text{a.e. in $\R_+$ for any $t\in \R_+$\,,}
\end{equation}
\begin{equation}\label{sottsass24}
v(\xi_n(t)^-,t)\,=\,v(\xi_n(t)^+,t)\,, \quad v_\xi(\xi_n(t)^-,t)
\,=\, b_n
v_\xi(\xi_n(t)^+,t)
\qquad (n\in\N, t\in [0,t_n))\,.
\end{equation}
It is natural to consider $v$ satisfying, beside \eqref{sottsass01}-\eqref{sottsass24}, the following initial-boundary conditions:
\begin{equation}\label{sottsass21}
v_{\xi}(0^+\!,t) \,=\,0 \quad\text{for all $t\in\R_+$}\,,\quad v(\xi,0)\,=\, v_0(\xi) \;\;\text{for a.e. $\xi\in\R_+$}
\end{equation}
\end{subequations}
with $v_0\in L^2(\R_+;\b)$, $0\le v_0\le 1$\,.

Super- and subsolutions of problem \eqref{sottsass} have properties \eqref{sottsass01}-\eqref{sottsass20} and satisfy \eqref{sottsass1}, the second relationship in \eqref{sottsass24} and \eqref{sottsass21} with proper inequalities. In particular, a stationary subsolution $\overline{p}$ of  \eqref{sottsass} satisfies
\begin{subequations}\label{ppp}
\begin{equation}\label{ppp0}
\overline{p}_n\in H^2(I_n) \,\textrm{ for all $n\in \N$}, \;
 \quad  \sum_{n\in \N}\int_{l_n}| \overline{p}_n''(\r)|^2\b(\r)\,d\r<\infty \,,
\end{equation}
\begin{equation}\label{ppp1}
\overline{p}'' + c\,\overline{p}' + f(\overline{p})\,\ge\,0  \quad\text{a.e. in $\R_+$\,,}
\end{equation}
\begin{equation}\label{ppp24}
\overline{p}'(0^+) \,\ge\,0 \,,\quad \overline{p}(\r_n^-)\,=\, \overline{p}(\r_n^+)\,, \quad \overline{p}'(\r_n^-) \,\le\, b_n\, \overline{p}'(\r_n^+) \,,
\end{equation}
\end{subequations}
whereas a stationary supersolution satisfies \eqref{ppp} with reverse inequalities.

The following comparison theorem will be proven.
\begin{theorem}\label{gecoxi}
Let $(H_1)$-$(H_2)$ be satisfied. Let $\underline{v}, \,\overline{v}$ be a sub-, respectively a supersolution of problem \eqref{sottsass}.  Then there holds  $\underline{v}(\xi,t) \le \overline{v}(\xi,t)$ for any $(\xi,t)\in Q_+$.
\end{theorem}
\begin{remark}\label{dacitare}
A  result analogous to Theorem \ref{decre}, whose formulation is left to the reader, follows from Theorem \ref{gecoxi}.
In particular, if $\underline{p}$ is a stationary subsolution of \eqref{sottsass}, the solution of \eqref{sottsass} with initial data $v(\cdot,0)=\underline{p}$ is nondecreasing in time and converges to a nontrivial weak stationary solution of \eqref{sottsass}. Arguing as in Remark \ref{uniso=1} shows that the only nontrivial weak stationary solution of \eqref{sottsass} is identically equal to $1$.
\end{remark}

To prove Theorem \ref{gecoxi} it is more convenient to work in the usual space $L^2(\R_+)$ than in the weighted space $L^2(\R_+;\b)$. Therefore, we make use of the unitary transformation
\begin{equation*}
U: L^2(\R_+;\b) \mapsto L^2(\R_+)\,,\qquad \tilde f \mapsto U\tilde f:= \sqrt{\b}\,\tilde f
\end{equation*}
(see \eqref{truni}). Let $\underline{v}$ be a subsolution of \eqref{sottsass}. Then the function $\underline{V}:= \sqrt{\b}\,\underline{v}$ belongs to $C(\overline{\R}_+;L^2(\R_+))\cap C^1(\R_+;L^2(\R_+))$,  $\underline{V}(\xi,t)=\sum_{n=1}^\infty \underline{V}_n(\xi,t)\chi_{_{J_n(t)}}(\xi)$, and the following holds:
\begin{subequations}\label{basta}
\begin{equation}\label{basta1}
  \underline{V}_n(\cdot,t)
  \in H^2(J_n(t))\;\textrm{ for all $n\in\N$}\,,\quad \sum_{n=1}^\infty\int_{J_n(t)}
   \Big(\big|\underline{V}_{n\xi\xi} (\cdot,t) \big|^2+\big|\underline{V}_{n\xi} (\cdot,t) \big|^2\Big)\,d\xi\,<\,\infty \quad\text{for any $t\in \R_+$\,,}
\end{equation}
\begin{equation}\label{sotto1}
\underline{V}_t(\cdot,t)\,\le\,  \underline{V}_{\xi\xi}(\cdot,t)\,+\, c\,\underline{V}_{\xi}(\cdot,t)\,+\,\sqrt\b f\left(\frac{\underline{V}(\cdot,t)}{\sqrt\b}\right)  \quad\text{a.e. in $\R_+$ for any $t\in \R_+$\,,}
\end{equation}
\begin{equation}
\underline{V}_{\xi}(0^+\!,t) \,\ge\,0 \;\;\text{for any $t\in \R_+$\,,}
\quad \underline{V}(\xi,0)\,\le\, \sqrt\b\,v_0(\xi) \;\;\text{for a.e. $\xi\in \R_+$}\,,
\end{equation}
\begin{equation}\label{solpro24}
\underline{V}_\xi(\xi_n(t)^-,t) \,\le\, \sqrt{b_n}\,\underline{V}_\xi(\xi_n(t)^+,t)  \quad (n\in\N, t\in \R_+)\,,
\end{equation}
\begin{equation}\label{solpro23}
\underline{V}(\xi_n(t)^+,t)\,=\, \sqrt{b_n}\,\underline{V}(\xi_n(t)^-,t)
\end{equation}
\end{subequations}
(see Subsubsection \ref{ss32}; observe that equality \eqref{solpro23} corresponds to the continuity of $\underline{v}(\cdot,t)$ at $\xi_n(t)$ in \eqref{sottsass24}). The function $\overline{V}:= \sqrt{\b}\,\overline{v}$, with $\overline{v}$ supersolution of \eqref{sottsass}, satisfies  \eqref{sotto1}-\eqref{solpro24} with reverse inequalities.

Observe that proving the ordering between $\underline{v}$ and $\overline{v}$ stated in Theorem \ref{gecoxi} is the same as proving it between $\underline{V}$ and $\overline{V}$, since the branching function $\b$ is constant in each interval $J_n(t)$.

For any $n\in\N$ and any fixed $\t\in(0, t_n]$, $t_n:=\frac{\r_n}{c}$\,, set
$$
K_n:=\{(\xi,t)\,|\, \xi\in(\xi_{n-1}(t),\xi_n(t)),\, t\in[0, t_n]\}\,, \quad K_{n,\t}:=K_n\cap\{(\xi,t)\,|\, \xi\in\R_+,\, t\in[0, \t]\}\,.
$$
The following lemma will be used.
\begin{lem}\label{backsol}
Let $a=a(\xi,t)\in L^\infty(Q_+)$ be given. Let  $\t\in(0, t_n]$ be fixed, and let $\eta=\eta(\xi)=\sum_{n=1}^\infty \eta_n(\xi)\chi_{_{J_n(\t)}}\!(\xi)$ with $\eta_n \in C_0^\infty(J_n(\t))$, $0\le\eta_n\le 1$ $(n\in\N)$. Then there exists $\varphi:\R_+\times[0,\t]\mapsto\R_+$\,, $\varphi(\xi,t)=\sum_{n=1}^\infty \varphi_n(\xi,t)\chi_{_{K_{n,\t}}}\!(\xi,t)$ such that for all $t\in(0,\t)$
\begin{subequations}\label{duffop}
\begin{equation}\label{duffop1}
\text{$\varphi_n (\cdot,t)\in L^2(\xi_{n-1}(t),\xi_{n+1}(t))$, \quad $\varphi_n(\cdot,t)|_{J_n(t)}\in H^2(J_n(t))$\,, \quad $\varphi_n(\cdot,t)|_{J_{n+1}(t)}\in H^2(J_{n+1}(t))$\,,}
\end{equation}
\begin{equation}\label{duffop3}
\varphi_{nt}(\cdot,t) \,=\, - \,\varphi_{n\xi\xi}(\cdot,t)+ c\,\varphi_{n\xi}(\cdot,t)-a(\cdot,t)\varphi_n(\cdot,t)
\quad\text{a.e. in $(\xi_{n-1}(t),\xi_{n+1}(t))$\,,}
\end{equation}
\begin{equation}\label{duffop2}
\varphi_{n\xi}(0^+\!,t)\,=\,0\,, \quad \varphi_n(\cdot,\t) =\eta_n\,,
\end{equation}
\begin{equation}\label{duffop4}
\varphi_{n+1}(\xi_n(t),t)= \sqrt{b_n}\,\varphi_n(\xi_n(t),t)\,,
\end{equation}
\begin{equation}\label{duffop5}
\varphi_{n\xi}(\xi_n(t),t)\,=\, \sqrt{b_n}\,(\varphi_{n+1})_\xi(\xi_n(t),t) \,.
\end{equation}
\end{subequations}
\end{lem}
\begin{proof}
Let $\phi$ be the solution of the backward problem
\begin{equation*}
\left\{
\begin{array}{ll}
\phi_t \,=\, A \phi -\tilde{a}\,\phi \quad\textrm{in}\,\,I_n\times (0,\t)
\\& \\
\phi = \tilde\eta
\quad \textrm{in}\,\, I_n\times \{\t\} \,,
\end{array}
\right.
\leqno{(AP)}
\end{equation*}
where $A$ is the operator characterized by Proposition \ref{colf} and $\tilde{a}=\tilde{a}(\r,t):=a(\r-ct,t)$, $\tilde\eta=\tilde\eta(\r):=\eta(\r-c\t)$. By standard semigroup theory there exists a unique nonnegative solution of problem $(AP)$ with $\phi\in C([0, \t];L^2( I_n))\cap C^1([0, \t);L^2( I_n))$, $\phi(\r,t)=\sum_{n=1}^\infty \phi_n(\r,t)\chi_{_{I_n}}(\r)$, $A\phi \in C([0, \t);L^2( I_n))$. Then by \eqref{duffop} for any $n\in\N$ there holds $\phi_n(\cdot,t)\in L^2(I_n)$ for any $t\in[0,\t]$, and  for all $t\in[0,\t)$ $\phi_n(\cdot,t)\in H^2(I_n)$\,,
\begin{equation*}
\phi_{nt}(\cdot,t) \,=\,-\, \phi_{n\r\r}(\cdot,t) \,+\,a(\cdot,t)\phi_n(\cdot,t)\quad\text{a.e. in $I_n$}\,,
\end{equation*}
\begin{equation*}
\phi_{n\r}(0^+\!,t)=0\,, \quad \phi(\r_n^+,t)= \sqrt{b_n}\,\phi(\r_n^-,t)\,,\quad\phi_{n\r}(\r_n^-,t)\,=\, \sqrt{b_n}\,\phi_{n\r}(\r_n^+,t)\,.
\end{equation*}
\smallskip
Setting $\varphi_n(\xi,t):=\phi_n(\xi+ct,t)$ $(\xi>0,\, 0\le t\le\t\le t_n; n\in\N)$ the claim follows.
\end{proof}

\smallskip

Now we can prove Theorem \ref{gecoxi}.

\smallskip

 \noindent{\em Proof of Theorem \ref{gecoxi}.}
$(i)$ Set $w:=\underline{V}-\overline{V}$, $w(\xi,t)=\sum_{n=1}^\infty w_n(\xi,t)\chi_{_{J_n(t)}}(\xi)$, and
$$
 \qquad a=a(\xi,t):=
\left\{
\begin{array}{ll}
\frac{f(\underline{v}(\xi,t))-f(\overline{v}(\xi,t)}{\underline{v}(\xi,t)-\overline{v}(\r,t)} \quad\text{if $w(\xi,t)\neq0$\,,}
\\& \\
0 \quad \text{otherwise}
\end{array}
\right.
$$
for all $(\xi,t)\in Q_+$ (observe that $a\in L^\infty(Q_+)$ by assumption $(H_2)$). Then by \eqref{basta} there holds $w(\cdot,0)\le0$ a.e. in $\R_+$, and for any $t>0$
\begin{subequations}\label{fine} and $n\in\N$
\begin{equation}\label{fine1}
w_{nt}(\cdot,t)-w_{n\xi\xi}(\cdot,t)- c\,w_{n\xi}(\cdot,t)- a(\cdot,t)w_n(\cdot,t)\,\le\, 0 \quad\text{a.e. in $J_n(t)$}\,,
\end{equation}
\begin{equation}\label{fine2}
w_{n+1}(\xi_n(t),t)\,=\,\sqrt{b_n}\,w_n(\xi_n(t),t)\,,
\end{equation}
\begin{equation}\label{solpro24}
w_{n\xi}(\xi_n(t),t)\,\le\, \sqrt{b_n} (w_{n+1})_\xi(\xi_n(t),t) \,,
\end{equation}
\begin{equation}\label{fine3}
w_{n\xi}(0^+\!,t) \,\ge\,0\,.
\end{equation}
\end{subequations}

Let $\t\in\R_+$ be fixed. For any $t\in[0,\t]$ there holds
\begin{eqnarray*}
\frac{d}{dt} \int_{\xi_{n-1}(t)}^{\xi_n(t)}w_n(\xi,t)\varphi_n(\xi,t)\,d\xi &=&
-\,c\;[w_n(\xi_n(t),t)\varphi_n(\xi_n(t),t) -w_n(\xi_{n-1}(t),t)\varphi_n(\xi_{n-1}(t),t)]\,+\\
&+&  \int_{\xi_{n-1}(t)}^{\xi_n(t)}[w_{nt}(\xi,t)\varphi_n(\xi,t)+ w_n(\xi,t)\varphi_{nt}(\xi,t)] \,d\xi \,.\nonumber
\end{eqnarray*}
Integrating the above equality on $[0,\t]$ gives
\begin{subequations}\label{calcoli}
\begin{eqnarray}\label{calcoli1}
&&  \iint_{K_{n,\t}} w_{nt}(\xi,t)\varphi_n(\xi,t)\,d\xi dt=\, - \iint_{K_{n,\t}} w_n(\xi,t)\varphi_{nt}(\xi,t)\,d\xi dt  \,+\\
&+&c\int_0^\t [ w_n(\xi_n(t),t)\varphi_n(\xi_n(t),t) -w_n(\xi_{n-1}(t),t)\varphi_n(\xi_{n-1}(t),t)]\,dt\,+ \nonumber\\
&+& \int_{\xi_{n-1}(\t)}^{\xi_n(\t)}w_n(\xi,\t)\eta_n(\xi)\,d\xi -\int_{\r_{n-1}}^{\r_n}w_n(\r,0)\varphi_n(\r,0)\,d\r \nonumber\,.
\end{eqnarray}
It is also easily seen that
\begin{eqnarray}\label{calcoli2}
&&\iint_{K_{n,\t}} w_{n\xi\xi}(\xi,t)\varphi_n(\xi,t)\,d\xi dt=\,  \iint_{K_{n,\t}}w_n(\xi,t)\varphi_{n\xi\xi}(\xi,t)\,d\xi dt\,+\\
&+&\int_0^\t [ w_{n\xi}(\xi_n(t),t)\varphi_n(\xi_n(t),t) -w_{n\xi}(\xi_{n-1}(t),t)\varphi_n(\xi_{n-1}(t),t)]\,dt\,- \nonumber\\
&-& \int_0^\t [w_n(\xi_n(t),t)\varphi_{n\xi}(\xi_n(t),t) -w_n(\xi_{n-1}(t),t)\varphi_{n\xi}(\xi_{n-1}(t),t)]\,dt\,,\nonumber
\end{eqnarray}
\begin{eqnarray}\label{calcoli3}
&&\iint_{K_{n,\t}} w_{n\xi}(\xi,t)\varphi_n(\xi,t)\,d\xi dt=\, - \iint_{K_{n,\t}}w_n(\xi,t)\varphi_{n\xi}(\xi,t)\,d\xi dt\,+\\
&+& \int_0^\t [w_n(\xi_n(t),t)\varphi_n(\xi_n(t),t) -w_n(\xi_{n-1}(t),t)\varphi_n(\xi_{n-1}(t),t)]\,dt\,. \nonumber
\end{eqnarray}
\end{subequations}
From \eqref{fine1} and \eqref{calcoli} we get
\begin{eqnarray*}
0&\ge& \iint_{K_{n,\t}} [w_{nt}-w_{n\xi\xi}- c\,w_{n\xi}- aw_n]\,\varphi_n\,d\xi dt \,=\\
&=& -\iint_{K_{n,\t}} [\varphi_{nt} + \varphi_{n\xi\xi}- c\,\varphi_{n\xi}+a\varphi_n]\,w_n\,d\xi dt\,+ \nonumber\\
&+& \int_{\xi_{n-1}(\t)}^{\xi_n(\t)}w_n(\xi,\t)\eta_n(\xi)\,d\xi -\int_{\r_{n-1}}^{\r_n}w_n(\r,0)\varphi_n(\r,0)\,d\r\,+ \nonumber\\
&-&\int_0^\t [ w_{n\xi}(\xi_n(t),t)\varphi_n(\xi_n(t),t) -w_{n\xi}(\xi_{n-1}(t),t)\varphi_n(\xi_{n-1}(t),t)]\,dt\,+ \nonumber\\
&+& \int_0^\t [w_n(\xi_n(t),t)\varphi_{n\xi}(\xi_n(t),t) -w_n(\xi_{n-1}(t),t)\varphi_{n\xi}(\xi_{n-1}(t),t)]\,dt\,. 
\end{eqnarray*}
Since $w(\cdot,0)\le0$ a.e. in $\R_+$, by \eqref{duffop3} the above inequality gives:
\begin{eqnarray*}
&& \int_{\xi_{n-1}(\t)}^{\xi_n(\t)}w_n(\xi,\t)\eta_n(\xi)\,d\xi \,\le \\
&\le& \int_0^\t [ w_{n\xi}(\xi_n(t),t)\varphi_n(\xi_n(t),t) - w_{n\xi}(\xi_{n-1}(t),t)\varphi_n(\xi_{n-1}(t),t)]\,dt\,- \nonumber\\
&-& \int_0^\t [  \varphi_{n\xi}(\xi_n(t),t)w_n(\xi_n(t),t) - \varphi_{n\xi}(\xi_{n-1}(t),t)w_n(\xi_{n-1}(t),t)]\,dt\,.\nonumber
\end{eqnarray*}

Let us now sum over $n\in\N$ the above family of inequalities; observe that the sum actually starts from $k\in\{1,\dots,n\}$ such that $t_{k-1}\le\t\le t_k$\,. Using \eqref{duffop2}-\eqref{duffop5} and \eqref{fine2}-\eqref{fine3} we get
\begin{eqnarray*}
&&\int_{\R_+}w(\xi,\t)\eta(\xi)\,d\xi \,\le \, -\int_0^\t [  w_{k\xi}(0^+,t)\varphi_k(0^+,t) - \varphi_{k\xi}(0^+,t)w_k(0^+,t) ]\,+\\
&+&\sum_{n=k}^\infty \int_0^\t [  w_{n\xi}(\xi_n(t),t)\varphi_n(\xi_n(t),t) - (w_{n+1})_\xi(\xi_n(t),t)\varphi_{n+1}(\xi_n(t),t)]\,dt\,-\\
&-&\sum_{n=k}^\infty \int_0^\t [ \varphi_{n\xi}(\xi_n(t),t)w_n(\xi_n(t),t) -(\varphi_{n+1})_\xi(\xi_n(t),t)w_{n+1}(\xi_n(t),t)]\,dt\,\le \\
&\le&\sum_{n=k}^\infty \int_0^\t [  w_{n\xi}(\xi_n(t),t)
-\sqrt{b_n}\, (w_{n+1})_\xi(\xi_n(t),t)]\varphi_n(\xi_n(t),t)\,dt\,- \\
&-&\sum_{n=k}^\infty \int_0^\t [ \varphi_{n\xi}(\xi_n(t),t)
-\sqrt{b_n}\,(\varphi_{n+1})_\xi(\xi_n(t),t)]w_n(\xi_n(t),t)\,dt\,\le\,0\,.
\end{eqnarray*}

By the arbitrariness of $\eta$, in the inequality
$$
\int_{\R_+}w(\xi,\t)\eta(\xi)\,d\xi \le 0
$$
we can choose $\eta=\eta_m$\,, $\lim_{m\to\infty}\eta_m= \chi_{\{w(\cdot,\t)>0\}}$ a.e. in $\xi\in\R_+$\,.
Then  we get
$$
\int_{\R_+}[w(\xi,\t)]_+\,d\xi  \le0\,,
$$
whence $w(\cdot,\t)\le0$ a.e. in $\R_+$\,. By the arbitrariness of $\t\in\R_+$ the conclusion follows.
\hfill $\square$


\section{Extinction versus propagation: Proofs}\label{expropro}
\setcounter{equation}{0}

Let us first prove Theorem \ref{t3}.

\noindent {\em Proof of Theorem \ref{t3}.}
 Let $\O_k\subseteq\mathcal{T}$ be an increasing family of subsets such that $O\in\O_k$, $\O_k\subset\subset \O_{k+1}$ for all $k\in\N$, and  $\bigcup_{k=1}^\infty \O_k=\mathcal{T}$. Consider the eigenvalue problem:
 \begin{equation*}
\left\{
\begin{array}{ll}
 \,  \D \varphi + [f'(0)+\mu]\varphi\,=\,0 &\textrm{in}\,\,\O_k \\& \\
\varphi \,=\,0& \textrm{on}\;\;\pa \O_k\setminus\{O\} \,.
\end{array}
\right.
\leqno{(EP_k)}
\end{equation*}
Denote by $\mu_{0,k}$ the first eigenvalue of $(EP_k)$ and by $\varphi_{0,k}=\varphi_{0,k}(x)>0$ the corresponding eigenfunction; also, let $\l_{0,k}$ denote the first eigenvalue of $-\D$ in $\O_k$ with Dirichlet zero boundary conditions on $\pa \O_k\setminus\{O\}$\,. Clearly,  there holds
\begin{equation}\label{stitreeter}
\mu_{0,k}\,=\,\lambda_{0,k} -f'(0)\,.
\end{equation}
Since by Proposition \ref{cokerdi}-$(ii)$ there holds $\lim_{k\to\infty}\l_{0,k}=E_0$\,, by \eqref{h5} and \eqref{stitreeter} there exists $k_0>0$ such that $\mu_{0,k_0}<0$\,.

For any $\e>0$ set
\begin{equation*}
w_{\e}(x)\,:=\,\left\{
\begin{array}{ll}
 \,  \e \varphi_{0,k_0}(x) &\textrm{in}\,\,\O_k\\& \\
 0& \textrm{otherwise.}
\end{array}
\right.
\end{equation*}
Since $\mu_{0,k_0}<0$\,, a standard calculation shows that for some $\e_0>0$ and
all $\e\in (0, \e_0)$ the function $w_{\e}$ is a stationary subsolution of equation \eqref{diffeq}.

On the other hand, by assumption $(H_2)$, the maximum principle and \eqref{tolex111} there holds $ u(\cdot,\bar t)>0 \; \mbox{in } \, \mathcal{T} $ for any fixed $\bar t>0$\,. Choosing $\e\in (0, \e_0)$ so small that $w_{\e}\le u(\cdot,\bar t)$ in $\mathcal{T}$,  by Theorem \ref{decre} we get
\begin{equation}\label{e34}
u_{\e}(\cdot,t)\;\le\, u(\cdot, t+\bar t)\;\le\, 1 \quad \text{in $\mathcal{T}$\; for all $t>0$\,,}
\end{equation}
where $u_{\e}$ denotes the solution of the problem
\begin{equation}\label{e32}
\left\{
\begin{array}{ll}
\, \displaystyle{ u_t  \,=\,  \D  u+f(u) }
 &\textrm{in}\,\,\mathcal{T}\times \R_+\\& \\
u \,=w_{\e}& \textrm{in}\;\;\mathcal{T}\times\{0\}\,.
\end{array}
\right.
\end{equation}

\noindent Moreover,

\noindent $(a)$ the function $t\mapsto u_{\e}(x,t)$ is nondecreasing in $\R_+$  for any $x\in \mathcal{T}$;

\noindent $(b)$ its pointwise limit $u_{\infty}(x):=\lim_{t\to \infty}\, u_{\e}(x,t)$ satisfies the equality
\begin{equation*}
\int_\mathcal{T} \{u_{\infty}\eta'' + f(u_{\infty})\eta\}\,d\m \,=\,0\quad\text{for all $\eta\in C^2_c(\mathcal{T})$;}
\end{equation*}

\noindent $(c)$  the convergence $u_{\e}(\cdot,t)\to u_{\infty}$ as $t\to \infty$ is  uniform on  compact subsets of $\mathcal{T}$ (see  Remark \ref{dinico}).

By Remark \ref{uniso=1} there holds $u_{\infty}\,\equiv\, 1$ in  $\mathcal{T}$, thus letting $t\to\infty$ in inequality \eqref{e34} the conclusion follows.
\hfill$\square$

\smallskip

The proofs of Theorems \ref{t30}-\ref{t2} rely on the construction of certain symmetric stationary supersolutions of equation \eqref{diffeq}. To this purpose, the ideas used in the proof of Theorem \ref{spehone} are expedient.
\begin{lem}\label{ln1}
Let  $\mathcal{T}$ be homogeneous, and let $\l\in(0,E_0)$. Then there exists $g=\sum_{n=1}^\infty g_n\chi_{_{I_n}}\in H^1(\R_+;\b)$, $0\le g\le1$  in $\R_+$\,, with the following properties:
\begin{subequations}\label{lisopra}
\begin{equation}\label{lisopra0}
g_n\in H^2(I_n) \,\textrm{ for all $n\in \N$}, \;
 \quad  \sum_{n\in \N}\int_{l_n}|g_n''(\r)|^2\b(\r)\,d\r<\infty \,,
\end{equation}
\begin{equation}\label{lisopra1}
 g_n''\,+\,\l g_n\,=\,0 \;\;\textrm{ in $I_n$ for all $n\in\N$\,,}
\end{equation}
\begin{equation}\label{lisopra2}
g'(0^+) \le0\,, \quad g_n'\,\le\,0 \;\textrm{ in $I_n$ for all $n\in\N$\,,}
\end{equation}
\begin{equation}\label{lisopra3}
g_{n+1}(\r_n)\,=\, g_n(\r_n) \,, \quad
g_n'(\r_n^-)\,\ge\, b\,g_{n+1}'(\r_n^+)\quad\text{for all $n\in\N$\,.}
\end{equation}
\end{subequations}
\end{lem}
\begin{proof}
Since $R\,:=\,\frac{b+1}{2\sqrt{b}}>1$ and by assumption $0<\sqrt{\l} r<\theta=\arccos{\frac 1R}$  (see \eqref{isp}), there holds $R\cos{(\sqrt{\l}r)}>1$, thus both roots \eqref{roqe} of equation \eqref{fqe} are positive. Observe that $\a_-<1$, and let $\a\in[\a_-,1)$ be fixed. Following \eqref{exqn} and  \eqref{truni}, we seek $ g_n$ of the following form:
\begin{equation}\label{rex}
 g_n(\r)\,:=\, \frac{\a^{n-1}}{b^{\frac{n}{2}}}\,\left\{\sqrt{b}\sin{[\sqrt{\l}(\r_n-\r)]} \,+\, \a\sin{[\sqrt{\l}(\r-\r_{n-1})]}\right\} \qquad(\r\in I_n)\,.
\end{equation}

It is easily seen that $g_n>0$  in $I_n$ (recall that by assumption $0<\sqrt{\l} r<\theta<\frac\pi 2$); moreover, equality \eqref{lisopra1} and the equality in \eqref{lisopra3} hold without restrictions on $\a$. In addition, there holds
\begin{equation}\label{lisopra22}
 g'(0^+)\, =\, g_1'(0^+)\, =\, \sqrt{\frac{\l}{b}}\,\left[\a- \sqrt{b}\cos{(\sqrt{\l}r)}\right]\,\le\,0\,,
\end{equation}
since
$$
\a\,<\,1\,<\,R\cos{(\sqrt{\l}r)}\,<\,\sqrt{b}\cos{(\sqrt{\l}r)}\,.
$$
This proves the first inequality in \eqref{lisopra2}, whence the second follows, since $ g_n'' =-\l  g_n\le0$ in $I_n$ for all $n\in\N$ (see \eqref{lisopra1}). Since $g_1(0)=\sin{(\sqrt{\l}r)}\le1$ and $g'\,\le\,0$ in $\R_+$\,, it also follows that $0\le g\le1$  in $\R_+$\,.
 Also inequality \eqref{lisopra3} is satisfied, since plainly by the choice of $\a\in[\a_-,1):$
\begin{eqnarray*}
g_n'(\r_n)\,-\, b\, g_{n+1}'(\r_n) &=&\sqrt{\l}\, \left(\frac{\a}{\sqrt{b}}\right)^{n-1}\left[-1 +\left( \sqrt{b}+\frac{1}{\sqrt{b}}\right)\cos{(\sqrt{\l}r)}\a-\a^2\right]
\\
&=& - \sqrt{\l}\,\left(\frac{\a}{\sqrt{b}}\right)^{n-1}  \left[\a^2 -2R\cos{(\sqrt{\l}r)}\a+ 1\right]\,\ge\,0\,. \nonumber
\end{eqnarray*}

Finally, from \eqref{lisopra1}, \eqref{rex} and the above remarks we get
\begin{equation}\label{liqinf}
0\,\le\, -\,\frac{g_n''}{\l}\,=\, g_n \,\le \,g_n(\r_{n-1}) \,= \, \left(\frac{\a}{\sqrt{b}}\right)^{n-1}\quad\text{ in $I_n$\,,}
\end{equation}
whence
$$
\frac{1}{\l^2} \sum_{n\in \N}\int_{l_n}(g_n'')^2\,\b\,d\r \,=\, \sum_{n\in \N}\int_{l_n}g_n^2\,\b\,d\r \,\le\,\sum_{n=1}^\infty\a^{2(n-1)}\,=\,\frac{1}{1-\a^2}\,.
$$
A similar calculation shows that
$$
\sum_{n\in \N}\int_{l_n} (g_n')^2\,\b\,d\r \,\le\,\left(1+\frac{\a}{\sqrt{b}}\right)^2\frac{\l}{1-\a^2}\,.
$$
Therefore, $g\in H^1(\R_+;\b)$ and \eqref{lisopra} hold true. Then the result follows.
\end{proof}

In the limiting case $\l=E_0$ the following holds:
\begin{lem}\label{lestali}
Let  $\mathcal{T}$ be homogeneous. Then there exists $h=\sum_{n=1}^\infty h_n\chi_{_{I_n}}$ with $h_n\in C^2(I_n)\cap C^1(\overline{I}_n)$ for all $n\in \N$, $0< h<\sqrt b$  in $\R_+$\,, such that
\begin{subequations}\label{ssopra}
\begin{equation}\label{ssopra1}
h_n''\,+\,E_0\, h_n\,=\,0 \;\;\textrm{ in $I_n$ for all $n\in\N$\,,}
\end{equation}
\begin{equation}\label{ssopra2}
h'(0^+) \le0\,, \quad  h_n'\,\le\,0 \;\textrm{ in $I_n$ for all $n\in\N$\,,}
\end{equation}
\begin{equation}\label{ssopra3}
h_{n+1}(\r_n)\,=\, \sqrt{b}\,h_n(\r_n)\,, \quad  h_n'(\r_n^-)\,=\, \sqrt{b}\,h_{n+1}'(\r_n^+)\quad\text{for all $n\in\N$\,.}
\end{equation}
\end{subequations}
\end{lem}
\begin{proof}
Since $\sqrt{E_0} r=\theta=\arccos{\frac 1R}$  (see \eqref{isp}), there holds $R\cos{(\sqrt{E_0}r)}=1$, thus $\a=1$ is the unique root of equation \eqref{fqe}. As in the proof of Lemma \ref{ln1}, this suggests to
seek a solution of \eqref{lestali} of the form $h=\sum_{n=1}^\infty h_n\chi_{_{I_n}}$\,, where
\begin{equation}\label{rexqn}
h_n(\r)\,:=\, \sqrt{b}\sin{\left[\sqrt{E_0}(\r_n-\r)\right]} \,+\, \sin{\left[\sqrt{E_0}(\r-\r_{n-1})\right]} \qquad(\r\in I_n)\,.
\end{equation}
Since $R>1$ there holds $\theta\in(0,\tfrac{\pi}{2})$, thus $\sin{[\sqrt{E_0}(\r_n-\r)]}>0$ if $\r\in[\r_{n-1},\r_n)$ and $\sin{[\sqrt{E_0}(\r-\r_{n-1})]}>0$ if $\r\in(\r_{n-1},\r_n]$. It follows that $h_n>0$ in the closed interval $\bar{I}_n$ for all $n\in\N$, thus $h>0$ in $\overline{\R}_+$\,.

Clearly, \eqref{ssopra1} is satisfied, and there holds
$$
h_{n+1}(\r_n) \,=\,\sqrt{b}\sin{\theta} \,=\,\sqrt{b}\,h_n(\r_n) \,,
$$
thus the first equality in \eqref{ssopra3} holds true. As for the second, from \eqref{rexqn} plainly we get
\begin{equation}\label{dexqn}
h'_n(\r)\,:=\,\sqrt{E_0} \,\left\{\,- \,\sqrt{b}\cos{\left[\sqrt{E_0}(\r_n-\r)\right]} \,+\, \cos{\left[\sqrt{E_0}(\r-\r_{n-1})\right]}\right\} \qquad(\r\in I_n)\,,
\end{equation}
whence
$$
   h_n'(\r_n)\,=\, \sqrt{E_0} \,\left[\,- \,\sqrt{b}\,+\, \cos{\theta}\right] \,,
 \quad
 h_{n+1}'(\r_n)\,=\, \sqrt{E_0} \,\left[\,- \,\sqrt{b} \cos{\theta}\,+\,1\right] \,.
$$
Then the second equality in  \eqref{ssopra3} holds by the very definition of $\theta$, since
\begin{eqnarray*}
 h_n'(\r_n^-)\,=\, \sqrt{b}\,h_{n+1}'(\r_n^+)& \Longleftrightarrow& - \,\sqrt{b}\,+\, \cos{\theta}\,=\,\,\sqrt{b} \,\left[\,- \,\sqrt{b} \cos{\theta}\,+\,1\right] \smallskip \\
& \Longleftrightarrow & \left(\sqrt{b} \,+\,\frac{1}{\sqrt{b}}\right) \cos{\theta}\,=\,2\\
& \Longleftrightarrow & \cos{\theta}\,=\,\tfrac1R\,.
\end{eqnarray*}
Finally,  \eqref{ssopra2} is satisfied since
$$
h'(0^+)\,=\,h_1'(0^+)\,=\, \sqrt{E_0} \, \left[\,- \,\sqrt{b} \cos{\theta}\,+\,1\right] \,=\, -\frac{b-1}{b+1}\,\sqrt{E_0}\,<\,0\,,
$$
and $h_n'' =-E_0\,  h_n\le0$ in $I_n$ for all $n\in\N$ (see \eqref{ssopra1}). From  \eqref{ssopra2} we obtain that $\|h\|_\infty=h(0)= \sqrt{b}\sin{\theta}<\sqrt b$ since $\theta\in(0,\tfrac{\pi}{2})$. Then the result follows.
 \end{proof}

Now we can prove Theorem \ref{t30}.

\smallskip

\noindent {\em Proof of Theorem \ref{t30}.}
Let $\l\in(f'(0),E_0)$ be fixed, and let $g=g(\r)$ be given by Lemma \ref{ln1}.
Set $\overline{z}:\overline{\R}_+\times\overline{\R}_+\mapsto [0,1]$, $\overline{z}(\r,t)\,:=\, g(\r)\,e^{-[\l-f'(0)]t}$. Plainly,  for any $t\in \R_+$ there holds
\begin{eqnarray*}
&&\overline{z}_t(\cdot,t) -  \overline{z}_{\r\r}(\cdot,t) - f(\overline{z}(\cdot,t))\,=\, -[\l-f'(0)]\,\overline{z}(\cdot,t)\,-\,g''\,e^{-[\l-f'(0)]t} - f(\overline{z}(\cdot,t))\, =\\
&=& - \left(g''+\l g\right)\,e^{-[\l-f'(0)]t} \,+\, f'(0)\,\overline{z}(\cdot,t)- f(\overline{z}(\cdot,t)) \,=\, f'(0)\,\overline{z}(\cdot,t)- f(\overline{z}(\cdot,t))\,\ge\,0
\end{eqnarray*}
a.e. in $\R_+$\,; here use of \eqref{lisopra1} and assumption $(H_4)$ has been made.

Let $u_0(x)\le g(\r(x))$ $(x\in\mathcal{T)}$. Then by Lemma \ref{ln1} and the above remarks the function
$\overline{u}(x,t):=\overline{z}(\r(x),t)$ $(x\in\mathcal{T}, t\in\overline{\R}_+)$ is a symmetric supersolution of problem \eqref{CN} (see Subsection \ref{symm}, whereas  $\underline{u}\equiv 0$ is clearly a subsolution. Therefore, by Theorem \ref{geco} the solution $u$ of \eqref{CN} satisfies $0\le u \le \bar{u}$ in $\mathcal{T}\times\R_+$\,. Since $\l>f'(0)$ and $\|g\|_{\infty}\le1$, we obtain that
$$
\lim_{t\to\infty}\|{u}(\cdot,t)\|_{\infty} \,\le\,\lim_{t\to\infty} e^{-[\l-f'(0)]t}\,=\,0\,.
$$
Hence the conclusion follows.
\hfill$\square$

\smallskip

The following proof is modeled after that of \cite[Theorem 5.2]{BPoT}, the function $\tilde{h}$ in \eqref{detih} below playing the role of a ground state.

\smallskip

\noindent {\em Proof of Theorem \ref{t2}.}  By assumption $(H_5)$ there exist $p_0\in(1,p)$ and $\s \in (0,1)$ such that
$$
f(u)\,\le\, E_0\, u^{p_0} \quad\text{ for any $u \in (0,\s)$}\,.
$$
Let $h$ be given by Lemma \ref{lestali}, and set
\begin{equation}\label{detih}
\tilde{h}(\r)\,:=\,k\,\frac{h(\r)}{\sqrt{\b(\r)}}\,, \quad\text{with $0\,<\,k\,<\,\frac {1} {\sqrt{b}} \,\left(\frac{\s}{1+\s}\right)^{\frac 1{p_0-1}}$}
 \qquad (\r\in\R_+)\,,
\end{equation}
Since $\|h\|_\infty< \sqrt{b}$ and $\b\ge1$ in $\R_+$\,, by \eqref{detih} there holds
$$
\|\tilde{h}\|_\infty\,\le\,k \sqrt{b}\,<\,\left(\frac{\s}{1+\s}\right)^{\frac 1{p_0-1}}<\,1\,.
$$
Moreover:
\smallskip

\noindent - by the first inequality in  \eqref{ssopra2} there holds $\tilde{h}'(0^+) \le0$;

\smallskip

\noindent - by the first equality in \eqref{ssopra3} and \eqref{detih} we get for all $n\in\N$:
\begin{equation*}
\tilde{h}(\r_n^-)\,=\, \tilde{h}_n(\r_n)\,=\,k\,\frac{h_n(\r_n)}{b^{\frac{n-1}{2}}}\,=\,k\,\frac{h_{n+1}(\r_n)}{b^{\frac{n}{2}}}
\,=\,\tilde{h}_{n+1}(\r_n)\,=\,\tilde{h}(\r_n^+)\,,
\end{equation*}
thus $\tilde{h}$ is continuous in $\overline{\R}_+$\,;

\smallskip

\noindent - by the second equality in \eqref{ssopra3} and \eqref{detih} there holds
\begin{equation*}
\tilde{h}'(\r_n^-)\,=\, \tilde{h}_n'(\r_n^-)\,=\,k\,\frac{h_n'(\r_n^-)}{b^{\frac{n-1}{2}}}\,=\,k\,\frac{h_{n+1}'(\r_n^+)}{b^{\frac n2-1}}
\,=\,\frac{b^{\frac{n}{2}}\tilde{h}_{n+1}'(\r_n^+)}{b^{\frac{n}{2}-1}} \,=\,b\,\tilde{h}_{n+1}'(\r_n^+)\,=\,b\,\tilde{h}'(\r_n^+)\,.
\end{equation*}

Now define $\overline{z}:\overline{\R}_+\times\overline{\R}_+\mapsto [0,1]$,
$$
\overline{z}(\r,t)\,:=\, \tilde{h}(\r)\,\z(t)\,e^{-E_0t}\,,
\quad\text{with
$ \zeta(t):=\left \{1-\|\tilde{h}\|_\infty^{p_0-1}\left[ 1- e^{-(p_0-1)E_0 t}\right]\right\}^{-\frac 1{p_0-1}}$\,}
\qquad (t\in\overline{\R}_+)\,.
$$
Plainly, for any $(\r,t)\in\overline{\R}_+\times\overline{\R}_+$ there holds
$$
0\,\le\,\overline{z}(\r,t)\, \le \,\tilde{h}(\rho)\,\zeta(t) \,\le\, \|\tilde{h}\|_\infty\,\|\z\|_{\infty}
\,=\,\frac{\|\tilde{h}\|_\infty}{\left[1-\|\tilde{h}\|_\infty^{p_0-1}\right]^{\frac 1{p_0-1}}}\,<\,\s \,,
$$
whence
\begin{equation}\label{disul}
f(\overline{z}(\r,t))\, \le \, E_0\, [\overline{z}(\r,t)]^{p_0}\,.
\end{equation}
Therefore,  for any $t\in \R_+$
\begin{eqnarray*}
&&\overline{z}_t(\cdot,t) -  \overline{z}_{\r\r}(\cdot,t) - f(\overline{z}(\cdot,t))\,\ge\,\overline{z}_t(\cdot,t) -  \overline{z}_{\r\r}(\cdot,t) - E_0\, [\overline{z}(\cdot,t)]^{p_0} \, =\smallskip\\
&=&-E_0\, \overline{z}(\cdot,t) + E_0\, \tilde{h}\, \|\tilde{h}\|_\infty^{p_0-1}\, e^{-p_0E_0 t}\z^{p_0}(t) -\tilde{h}''\z(t)\,e^{-E_0t}
-E_0\, \tilde{h}^{p_0}\, e^{-p_0E_0 t}\z^{p_0}(t)
\,= \smallskip\\
&=& - \left(\tilde{h}''+E_0 \tilde{h}\right)\z(t)\,e^{-E_0t} \,+\,
E_0\, \tilde{h}\, \left(\|\tilde{h}\|_\infty^{p_0-1}-\tilde{h}^{p_0-1}\right)\, e^{-p_0E_0 t}\z^{p_0}(t)\,\ge\,0
\end{eqnarray*}
a.e. in $\R_+$\,; here use of \eqref{ssopra1} and  \eqref{disul} has been made.

Let $u_0(x)\le \tilde{h}(\r(x))$ $(x\in\mathcal{T)}$. Then by the above considerations the function
$\overline{u}(x,t):=\overline{z}(\r(x),t)$ $(x\in\mathcal{T}, t\in\overline{\R}_+)$ is a symmetric bounded supersolution of problem \eqref{CN}, whereas  $\underline{u}\equiv 0$ is clearly a subsolution. Therefore, by Theorem \ref{comparison} the solution $u$ of \eqref{CN} satisfies $0\le u \le \bar{u}$ in $\mathcal{T}\times\R_+$\,. Since
$$
\sup_{x\in\mathcal{T}}\, \overline{u}(x,t)\,\le\,
 \|\tilde{h}\|_\infty\,\|\z\|_{\infty}\,e^{-E_0 t}
\,=\,\frac{\|\tilde{h}\|_\infty}{\left[1-\|\tilde{h}\|_\infty^{p_0-1}\right]^{\frac 1{p_0-1}}}\,e^{-E_0 t}\, \to\, 0\quad\textrm{as}\,\,t\to \infty,
$$
the conclusion follows.
\hfill $\square$


\section{Asymptotic speed of propagation: Proofs}\label{aspropro}
\setcounter{equation}{0}

\subsection{Proof of Theorem \ref{t1}}
Two preliminary lemmata are needed.
\begin{lem}\label{lesta0}
Let  the assumptions of Theorem \ref{t1} be satisfied. Then  there exists a function $m:=\sum_{n=1}^\infty m_n\chi_{_{I_n}}$, $0\le m\le1$ in $\R_+$\,, with the following properties:
\begin{subequations}\label{apro}
\begin{equation}\label{apro0}
m_n\in C^2(I_n)\cap C^1(\overline{I}_n)\;\textrm{ for all $n\in\N$}\,,
 \end{equation}
\begin{equation}\label{apro1}
 m_n''+c m_n'+f(m_n)\,\le\,0 \quad\textrm{in $I_n$\,, \,for all $n\in\N$ and $c>\hat c$}
\end{equation}
with $\hat c$ given by \eqref{decabis},
\begin{equation}\label{apro2}
m'(0^+) \le0\,,\quad m_n'\,\le\,0 \;\textrm{in $I_n$ for all $n\in\N$\,,}
\end{equation}
\begin{equation}\label{apro3}
m_{n+1}(\r_n)\,=\, m_n(\r_n) \,, \quad  m_n'(\r_n^-)\,\ge\, b_n\,m_{n+1}'(\r_n^+)\quad\text{for all $n\in\N$\,,}
\end{equation}
\begin{equation}\label{apro4}
\lim_{\r\to\infty} m(\r)=0\,.
\end{equation}
\end{subequations}
\end{lem}
\begin{remark}\label{loi}
By analogy with the case of \eqref{sproge}, the function $m$ given by Lemma \ref{lesta0} can be regarded as a supersolution of the {\em impulsive differential problem:}
\begin{subequations}\label{prosub}
\begin{equation}\label{pro21}
 q''+c q'+f(q)\,=\,0 \quad\textrm{a.e. in}\,\,\R_+\,,
\end{equation}
\begin{equation}
q'(0^+) =0\,,
\end{equation}
\begin{equation}\label{pro23}
q(\r_n^+)\,=\,q(\r_n^-)\,, \quad  q'(\r_n^-)\,=\, b_n\,q'(\r_n^+)\quad\text{for all $n\in\N$}
\end{equation}
\end{subequations}
(see \cite{LBS, RT1, RT2}). Subsolutions of \eqref{prosub} are similarly defined.
\end{remark}

\noindent {\em Proof of Lemma \ref{lesta0}.} Set $m_1:=1$ in the interval $I_1$, and for any $n\in\N$, $n\ge2$
\begin{subequations}\label{zx8}
\begin{equation}\label{zx8a}
m_n(\r)\,:=\,\frac{1}{\sqrt{b_0\dots b_{n-1}}}\left(\frac{1}{\sqrt{b_0\dots b_{n-2}}} \,-\, \kappa_n\,(\r-\r_{n-1})\right)  \qquad\, (\r\in I_n)
\end{equation}
with
\begin{equation}\label{zx8b}
 \kappa_n\,:=\, \frac{\frac{1}{\sqrt{b_0\dots b_{n-2}}}-\frac{1}{\sqrt{b_0\dots b_n}}}{\r_n-\r_{n-1}}\,>\,0 \,.
\end{equation}
\end{subequations}
Clearly, by definition $0<m\le1$ in $\R_+$ and \eqref{apro0} holds true. Moreover, for all $n\in\N$ the equality in \eqref{apro3} is satisfied. In fact,
$$
m_{n+1}(\r_n) \,=\,\frac{1}{\sqrt{b_0\dots b_n}}\,\frac{1}{\sqrt{b_0\dots b_{n-1}}}\,=\, \frac{1}{\sqrt{b_0\dots b_{n-1}}}\,\frac{1}{\sqrt{b_0\dots b_n}}\,=\,\,m_n(\r_n)\,.
$$

Since $m_1\equiv 1$ in $I_1$ by definition, there holds $m'(0^+)=m_1'(0^+) =0$\,, thus the first inequality in \eqref{apro2} holds true. The second is also satisfied, since $m_n'= -\, \frac{\kappa_n}{\sqrt{b_0\dots b_{n-1}}}$ in $I_n$ $(n\in\N)$. In addition, since
 $$
 m_n'(\r_n^-) \,=\,-\,\frac{\kappa_n}{\sqrt{b_0\dots b_{n-1}}}\,, \qquad  m_{n+1}'(\r_n^+) \,=\,-\,\frac{\kappa_{n+1}}{\sqrt{b_0\dots b_n}}\,,
 $$
the inequality in \eqref{apro3} follows if
$$
-\kappa_n\,\ge\, -\,\sqrt{b_n}\,\kappa_{n+1} \quad \Longleftrightarrow \quad \sqrt{b_n}\, \frac{\sqrt{b_n b_{n+1}}-1}{\sqrt{b_0\dots b_{n+1}}(\r_{n+1}-\r_n)} \;\ge\; \frac{\sqrt{b_{n-1}b_n}-1}{\sqrt{b_0\dots b_n}(\r_n-\r_{n-1})}\,.
$$
By assumption $(H_0)$-$(ii)$ there holds $\frac{1}{\r_{n+1}-\r_n} \;\ge \frac{1}{\r_n-\r_{n-1}}$\,, thus the latter inequality is satisfied if
$$
\sqrt{b_n}\,\left(\sqrt{b_n b_{n+1}}-1\right)\;\ge\;\sqrt{b_{n+1}}\,\left(\sqrt{b_{n-1}b_n}-1\right)\;,
$$
namely
$$
\sqrt{b_nb_{n+1}}\,\left(\sqrt{b_n}-\sqrt{b_{n-1}}\right)\;\ge\; \sqrt{b_n}- \sqrt{b_{n+1}} \,.
$$
The above inequality is clearly true, since by $(H_0)$ the sequence $\{b_n\}$ is nondecreasing.
Hence the inequality in \eqref{apro3} holds true.

\smallskip

Further, let us show that inequality \eqref{apro1} holds. This is obviously tue in $I_1$\,. For all $n\ge2$, since $m_n'=
-\frac{\kappa_n}{\sqrt{b_0\dots b_{n-1}}}$\,, $m_n''=0$, this amounts to prove that for all $n\in\N$ and $c>\hat c$
\begin{equation}\label{qaz1}
 f(m_n(\r)) \,\le \, c\; \frac{\kappa_n}{\sqrt{b_0\dots b_{n-1}}} \qquad \qquad (\r\in I_n)\,.
\end{equation}
Since $f(m_n(\r))\le Mm_n(\r)$ by assumption $(H_3)$, and
$$
m_n(\r)\,\le\, m_n(\r_{n-1})\,=\,\frac{1}{\sqrt{b_0\dots b_{n-1}}}\frac{1}{\sqrt{b_0\dots b_{n-2}}}\quad\text{for any $\r\in I_n$\,,}
$$
inequality \eqref{qaz1} holds if for all $n\in\N$
\begin{equation}\label{qaz2}
\frac{M}{\sqrt{b_0\dots b_{n-2}}}\, \le\, c\,\frac{\sqrt{b_{n-1} b_n}-1}{\sqrt{b_0\dots b_n}(\r_n-\r_{n-1})} \qquad \Longleftrightarrow \qquad
c\,\ge\, M(\r_n-\r_{n-1})\,\frac{\sqrt{b_{n-1} b_n}}{\sqrt{b_{n-1} b_n}-1}\,.
\end{equation}
If $c>\hat c$, it follows from $(H_0)$-$(ii)$ that the latter inequality holds for all $n\in\N$. This proves \eqref{qaz1}, whence\eqref{apro1} follows.

\smallskip

Concerning \eqref{apro4}, by \eqref{zx8} there holds
\begin{equation}\label{zx9}
m(\r)\,\le\, \chi_{_{I_1}}(\r)  + \sum_{n=2}^\infty \frac{1}{\sqrt{b_0\dots b_{n-1}}\sqrt{b_0\dots b_{n-2}}} \,\chi_{_{I_n}}(\r)\,=:\, F(\r) \qquad (\r\in\R_+)\,.
\end{equation}
Therefore,
\begin{equation}\label{zx10}
\lim_{\r\to\infty}m(\r)\,\le\, \lim_{\r\to\infty} F(\r)\,=\, \lim_{n\to\infty}\frac{1}{\sqrt{b_0\dots b_{n-1}}\sqrt{b_0\dots b_{n-2}}}
\,\le\, 2\sqrt2\, \lim_{n\to\infty}\frac{1}{2^n}\,=\,0 \,.
\end{equation}
This completes the proof.
\hfill $\square$

\begin{lem}\label{lesopra}
Let  the assumptions of Theorem \ref{t1} be satisfied. Let $c>\hat c$, with $\hat c$ given by \eqref{decabis}, and let $m$ be given by Lemma \ref{lesta0}. Let $u_0(x) \le m(\r(x))$ for a.e. $x\in\mathcal{T}$. Then the function
\begin{equation}\label{defiv}
U\equiv U_c:\mathcal{T}\times \R_+\mapsto[0,1]\,, \qquad
U_c(x,t)\,:=\,\left\{
\begin{array}{ll}
m(\r(x)-ct)\quad\text{if $\r(x)>ct$\,,}\\& \\
1\qquad\qquad\qquad\text{if $\r(x)\le ct$}
\end{array}
\right.
\end{equation}
 is a weak bounded supersolution of problem \eqref{CN}.
\end{lem}

\begin{proof} It is easily seen that $U\in C(\overline{\R}_+; H_{loc}^1(\mathcal T))\cap L^\infty(\mathcal{T}\times \R_+),$ since $m_n\in C^2(I_n)\cap C^1(\overline{I}_n)$  for any $n\in\N$ and $0\le m\le 1$ in $\R_+$\,. Concerning \eqref{e1000}, observe that 
\begin{subequations}\label{000}
\begin{equation}\label{0022}
U_t(x,t) \,=\, 0\,=\,f(1)\,=\, U_{xx}(x,t) +f(U(x,t)) 
\end{equation}
if $\r(x)\le ct$, whereas by  \eqref{apro1}
\begin{equation}\label{0021}
-c\,m'(\r(x)-ct) \,\ge\,
 m''(\r(x)-ct) \,+\,f(m(\r(x)-ct))
\end{equation}
 for a.e. $x\in\mathcal{T}$ and any $t\in\R_+$ such that $\r(x)>ct$. Moreover, by assumption there holds 
\begin{equation}\label{0023}
 u_0(x) \le m(\r(x))\le U(x,0) \quad\text{ for a.e. $x\in\mathcal{T}$}\,.
\end{equation}
\end{subequations}
  From \eqref{000} inequality \eqref{e1000} immediately follows.
 
Concerning \eqref{solprop23}, observe that the first inequality is satisfied, since $\r(O)=0$ and by \eqref{defiv} there holds $U(\cdot,t)=1$ in the interval $(0,ct)$ $(t\in\R_+)$. As for the second, it is trivially satisfied at any vertex of the $n$-th generation when $\r_n\le ct$\,. It remains to prove that 
\begin{equation}\label{mmm}
m'(\r_n^--ct)\,\ge\, b_n\,m'(\r_n^+-ct)\quad\text{for all $n\in\N$ and $t\in\left(0,\frac{\r_n}{c}\right)$\,.}
\end{equation}

Let $n\in\N$ be fixed. Let $t=t_k:=\frac{\r_n-\r_k}{c}$ $(k=1,\dots,n-1)$, thus $\r_n-ct_k=\r_k$ and \eqref{mmm} reads 
$$
m'(\r_k^-)\,\ge\, b_n\,m'(\r_k^+)\qquad (k=1,\dots,n-1)\,.
$$
 Since $b_n\ge b_k$ for all $k=1,\dots,n-1$ and $m'\le0$ in $\R_+$ (see $(H_0)$-$(ii)$ and \eqref{apro2}), from  \eqref{apro3} we get
$$
 m'(\r_k^-)\,=\, m_k'(\r_k^-)\,\ge\, b_k\,m_{k+1}'(\r_k^+)\,\ge\, b_n\,m_{k+1}'(\r_k^+)\,=\, b_n\,m'(\r_k^+)\qquad (k=1,\dots,n-1)\,,
$$
thus \eqref{mmm} follows when $t\in\{t_1,\dots,t_{n-1}\} \subseteq\left(0,\frac{\r_n}{c}\right)$\,.

Further observe that by \eqref{zx8a} and the equality in \eqref{apro3} there holds $m\in C(\overline{\R}_+)$, thus in particular
\begin{equation*}
m(\r_n^--ct)\,=\,m(\r_n^+-ct) \quad\text{for all $n\in\N$, $t\in\left(0,\frac{\r_n}{c}\right)$\,.}
\end{equation*}
For any $t\in\left(0,\frac{\r_n}{c}\right)\setminus\{t_1,\dots,t_{n-1}\}$ the above equality can be differentiated in the classical sense since $m_n\in C^1(\overline{I}_n)$, thus $m_n(\cdot-ct)\in C^1(\overline{J(t)}_n)$. Then for any $t\in\left(0,\frac{\r_n}{c}\right)\setminus\{t_1,\dots,t_{n-1}\}$ we obtain
\begin{equation*}
-c\,m'(\r_n^--ct)\,=\, -c\,m'(\r_n^+-ct)\quad \Rightarrow\quad m'(\r_n^--ct)\,=\, m'(\r_n^+-ct)\,\ge\, b_n\,m'(\r_n^+-ct)\,,
\end{equation*}
since $b_n\ge2$ and $m'\le0$ in $\R_+$\,. Then by the arbitrariness of $n\in\N$ the second inequality in \eqref{solprop23} follows. This completes the proof.
\end{proof}

\smallskip

Now we can prove Theorem \ref{t1}.
\smallskip

\noindent {\em Proof of Theorem \ref{t1}.}
Let $\bar{c}>\hat{c}$ be fixed. By Lemma \ref{lesopra} $U_{\bar{c}}$ is a weak  bounded  supersolution of problem \eqref{CN}, thus by Theorem \ref{comparison} (see also Remark \ref{wefo}) there holds $0\le u\le U_{\bar{c}}$ in $\mathcal{T}\times \R_+$. In particular, for any fixed $y\in\mathcal{T}$ we get
\begin{equation}\label{nnn}
u(x,t) \,\le\, m(\r(x)-\bar{c}t) \quad\text{for all $(x,t)\in\mathcal{T}\times \R_+$ with $\r(x)\ge\r(y)+ \bar{c}t$\,.}
\end{equation}

Let $\e\in(0,1)$ be arbitrarily fixed. Since the inequality $\r(x)\ge\r(y)+ \bar{c}t$ implies both $\r(x)>\r(y)+(1-\e)\,\bar{c}\,t$\, and $\r(x)-(1-\e)\,\bar{c}\,t\ge \r(y)+\e \bar{c}t\,,$ from \eqref{nnn} and the nonincreasing character of $m$ (see \eqref{apro2}) we obtain that
\begin{equation*}
u(x,(1-\e)t) \,\le\, m(\r(x)-(1-\e)\,\bar{c}\,t) \,\le\, m(\r(y)+\e\, \bar{c}\,t)
\quad\text{for all $(x,t)\in\mathcal{T}\times \R_+$ with $\r(x)\ge\r(y)+ \bar{c}t$\,.}
\end{equation*}
Therefore, for any $y\in\mathcal{T}$, $\e\in(0,1)$ and $t\in \R_+$ there holds
\begin{equation*}
\sup_{\r(x)\ge\r(y)+ \bar{c}t}
u(x,(1-\e)t) \,\le\,m(\r(y)+\e\, \bar{c}\,t)\,.
\end{equation*}
Set $\bar{c}_\e:=\frac{\bar{c}}{1-\e}\,,$ $s:=(1-\e)t$\,. Then the above inequality reads
\begin{equation}\label{bbb}
\sup_{\r(x)\ge\r(y)+ \bar{c}_\e s}
u(x,s) \,\le\,m\left(\r(y)+\e \,\bar{c}_\e\, s\right)\,.
\end{equation}

Letting $s\to\infty$ in \eqref{bbb} and using \eqref{apro4} gives \eqref{e89} with $c=\bar{c}_\e$\ for all  $\e\in(0,1)$. When $\e$ varies over $(0,1)$, $\bar{c}_\e$ varies over $(\bar{c},\infty)$. Therefore, we proved that for any fixed $y\in\mathcal{T}$ and  $\bar{c}>\hat{c}$ equality \eqref{e89} holds in the half-line $(\bar{c},\infty)$. Then by the arbitrariness of $\bar{c}$ the result follows.
\hfill$\square$


\subsection{Proof of Theorem \ref{t4}.} \label{teo25}

 Let us first prove some preliminary results.
\begin{lem}\label{lestasco}
Let  the assumptions of Theorem \ref{t4} be satisfied. Let $c\in(0,\check{c})$, $\m\in(-E_0,0)$,  and let $N\in\N$ be fixed.
Then there exists a function $\psi:=\sum_{n=1}^N \psi_n\chi_{_{I_n}}$, $0\le \psi\le1$ in $\R_+$, such that
\begin{subequations}\label{apropos}
\begin{equation}\label{apropos0}
\psi_n\in C^2(I_n)\cap C^1(\overline{I}_n)\;\;\textrm{ for all $n=1,\dots,N$},
 \end{equation}
\begin{equation}\label{apropos1}
 \psi_n''+c \psi_n'+[f'(0)+\m]\psi_n\,=\,0 \quad\textrm{in $I_n$ for all $n=1,\dots,N$,}
\end{equation}
\begin{equation}\label{apropos2}
\psi'(\r)\,<\,0 \;\quad\textrm{for any $\r\in(0,\r_N)$\,,}
\end{equation}
\begin{equation}\label{apropos3}
\psi_{n+1}(\r_n)\,=\, \psi_n(\r_n) \,, \quad  \psi_n'(\r_n^-)\,=\, b_n\,\psi_{n+1}'(\r_n^+) \quad\text{for all $n=1,\dots,N-1$\,,}
\end{equation}
\begin{equation}\label{apropos4}
\psi'(0^+) =0\,, \quad \psi(\r_N) =0\,.
\end{equation}
\end{subequations}
\end{lem}
\begin{proof}  Set $\o^2:=f'(0)+\m-\frac{c^2}{4}>0$ (see \eqref{altrecdz}), and
\begin{equation}\label{depsi}
\psi_n(\r)\,:=\, e^{-\frac{c\r}{2}}\big\{A_n\sin{\left[\o(\r-\r_{n-1})\right]\,+\, B_n\sin{\left[\o(\r_n-\r)\right]}}\big\}
 \qquad(\r\in I_n)
\end{equation}
with $A_n>0, \,B_n>0$ to be fixed $(n=1,\dots,N)$. Clearly, \eqref{apropos0}-\eqref{apropos1} are satisfied. Since by assumption $\r_1<\frac{\pi}{2\sqrt{f'(0)}}$\,, there holds $\o(\r_n-\r_{n-1})<\sqrt{f'(0)}\r_1<\frac\pi2$\,, thus $\psi>0$ in $\bigcup_{n=1}^NI_n$\,. Imposing the condition $\psi_1'(0^+) =0$ plainly gives
\begin{subequations}\label{detab}
\begin{equation}\label{detab1}
B_1\,=\, \frac{2\o}{c\,\sin{(\o\r_1)}+2\o\cos{(\o\r_1)}}\;A_1 \,.
\end{equation}
Since $\psi>0$ in $\bigcup_{n=1}^NI_n$, from the equalities $\psi_1'(0^+) =0$ and
\begin{equation}\label{negaderi}
\psi'(\r)\,=\,-\,[f'(0)+\m]\int_0^\r e^{-c(\r-\s)}\psi(\s)\,d\s \qquad (\r\in(0,\r_N))
\end{equation}
we get inequality \eqref{apropos2}
(observe that $f'(0)+\m> f'(0)-E_0>0$).

An elementary calculation shows that equalities \eqref{apropos3} read
\begin{equation}\label{detab2}
B_{n+1}\,=\, \frac{\sin{[\o(\r_n-\r_{n-1})]}}{\sin{[\o(\r_{n+1}-\r_n)]}}\;A_n\,,
\end{equation}
respectively
\begin{equation}\label{detab3}
A_{n+1} \,=\,\cos{[\o(\r_{n+1}-\r_n)]}\, B_{n+1}\,+
\,\frac{1}{b_n }\,\big\{\cos{[\o(\r_n-\r_{n-1})]}A_n
- B_n\big\}
\end{equation}
\end{subequations}
for any $n=1,\dots, N-1$. From \eqref{detab} all coefficients $B_1$ and $A_n, B_n$ with $n=2,\dots, N$ can be determined
depending on $A_1$\,. Finally, $A_1$ is determined imposing the condition $\psi_N(\r_N) =0$. Then the result follows.
\end{proof}

\begin{lem}\label{lestatilde}
Let  the assumptions of Theorem \ref{t4} be satisfied, and let $\psi$ be given by Lemma \ref{lestasco}. Then  for any $\e>0$ sufficiently small the function $\underline{\psi_{\e}}\!:\,\R_+\mapsto[0,1]$\,,
\begin{equation}\label{defipsie}
\underline{\psi_\e}(\r)
\,:=\,\left\{
\begin{array}{ll}
 \,  \e\, \psi(\r)&\textrm{if}\,\,\r\in(0,\r_N),
 \\& \\
 0& \textrm{otherwise,}
\end{array}
\right.
\end{equation}
is a stationary subsolution of problem  \eqref{sottsass}. Moreover, there holds $\underline{\psi_\e}\le0$ in $\R_+$\,.
\end{lem}
\begin{proof}
Let us check that \eqref{ppp1}-\eqref{ppp24} are satisfied. Set $\underline{\psi}\equiv\underline{\psi_{\e}}$ for simplicity.
From \eqref{apropos1} by a suitable choice of $\e$ we get
\begin{equation}\label{okm}
 \underline{\psi}''+c\,\underline{\psi} '+f(\underline{\psi} )\,=\, f(\e\psi)-[f'(0)+\m](\e\psi)\,=\,\e\psi \,[\,|\m|\,+\, o(\e)]\,> \, 0\quad\text{in $(0,\r_N)$}\,.
\end{equation}
To complete the proof, it suffices to observe that
$$
\underline{\psi}(\r_N^-)=\e\psi(\r_N)=0=\underline{\psi}(\r_N^+)\,, \quad
\underline{\psi}'(\r_N^-)=\e\psi'(\r_N)<0=\underline{\psi}'(\r_N^+)
$$
(see \eqref{negaderi}). Hence the result follows.
\end{proof}

Now we can prove  Theorem \ref{t4}.

\smallskip

\noindent {\em Proof of Theorem \ref{t4}.}
Let $\underline{\psi_{\e}}$ be given by Lemma \ref{lestatilde}, and let $\e>0$ be so small that $\underline{\psi_{\e}}(\r(x))\le u_0(x)$ for any $x \in  \mathcal{T}$. Then by Theorem \ref{geco} there holds
\begin{equation}\label{defiv1}
w(\r(x),t)\le u(x,t) \;\;\text{for any $x \in  \mathcal{T}$ and $t\in\R_+$}\,,
\end{equation}
where $u$ solves problem \eqref{CN} with initial data $u_0$ and $w$ solves problem \eqref{sottsass} with initial data $\underline{\psi_{\e}}$. Since $\underline{\psi_{\e}}'\le0$ (see Lemma \ref{lestatilde}), by a standard argument there holds $w_\r(\cdot,t)\le0$ in $\R_+$ for any $t\in\R_+$\,.

Set $v(\xi,t)\,:=\,w(\xi+ct,t)$\,. It is immediately seen that $v$ is a supersolution of problem \eqref{sottsass} with initial data $\underline{\psi_{\e}}$ (in this connection, observe that $v_\xi(0^+,t)=w_\r(ct,t)\le0$). Hence by Theorem \ref{gecoxi} there holds $v\ge z$ in $Q_+$\,, $z$ being the solution of problem \eqref{sottsass} with initial data $\underline{\psi_{\e}}$. On the other hand, since $\underline{\psi_{\e}}$ is a stationary subsolution of the same problem, by Remark \ref{dacitare} there holds
\begin{equation*}
\lim_{t\to\infty} z(\r(y),t)  \,=\,1 \quad\text{for any fixed $y \in  \mathcal{T}$}\,.
\end{equation*}

It follows that
\begin{equation}\label{defiv2}
\lim_{t\to\infty} v(\r(y),t) \,=\, \lim_{t\to\infty} w(\r(y)+ct,t) \,=\,1 \quad\text{for any fixed $y \in  \mathcal{T}$}.
\end{equation}
From \eqref{defiv1}-\eqref{defiv2} we obtain
\begin{equation*}
\lim_{t\to \infty} \inf_{\r(x)<\r(y) +c t}u(x,t)\,\ge\, \lim_{t\to \infty} \inf_{\r(x)<\r(y) +c t}w(\r(x),t)\,=\, \lim_{t\to \infty} w(\r(y) +c t,t)\,=\,1\,.
\end{equation*}
Then the conclusion follows.
\hfill$\square$


\section*{Appendix}\label{ar}
\renewcommand{\thesection}{\Alph{section}}\setcounter{section}{1}\setcounter{subsection}{0}\setcounter{equation}{0}


\subsection{Metric graphs}\label{megra}

While referring the reader to \cite{BK, Mu} for general definitions and results concerning metric graphs, in this subsection we recall some general concepts for convenience of the reader.

\subsubsection{General facts}\label{gefagra}

Like combinatorial graphs  (e.g., see \cite{Mu}), a metric graph consists of a countable set of vertices and of a countable set of edges. However, in contrast to combinatorial graphs, the edges are regarded as intervals glued together at the vertices:

\begin{definition}\label{meg} A {\em metric graph}  is a quadruple $(E, V, i, j)$, where:

\noindent $(a)$ $E$  is a countable family of open intervals $I_e\equiv(0,l_e)$, called {\em edges} of {\em length} $l_e\in(0,\infty]$;

\noindent $(b)$ $V$ is a countable set of points, called {\em vertices};

\noindent $(c)$ the maps $ i:E\mapsto V$ and $j:\{e\in E\,|\, l_e<\infty\} \mapsto V$ define the {\em initial point} and the {\em final point} of an edge. Both the initial  and the final point are {\em endpoints} of an edge.

\noindent A metric graph is {\em finite} if both $E$ and $V$ are finite.  A finite metric graph with $l_e<\infty$ for all $e\in E$ is called {\em compact}.
\end{definition}
We always assume for simplicity that $i(e)\neq j(e)$ (namely, no {\em loops} exist), and $l_e<\infty$ for all $e\in E$.
The following notations will be used:
$$
\mathcal{G}_e :=\{e\}\cup  I_e\,, \quad \mathcal{G}:= (\bigcup_{e\in E}\mathcal{G}_e)\cup V\,, \quad
 \overline{\mathcal{G}}_e:=\mathcal{G}_e \cup\{i(e),j(e)\}\,.
$$
By abuse of language, we currently speak of the metric graph $\mathcal{G}$.

\smallskip

We shall write $e\ni v$ or $v\in e$, if either $i(e)=v$ or $j(e)=v$ $(e\in E, v\in V)$. The notation $e\equiv(u,v)$ means that $i(e)=u$ and $j(e)=v$; therefore, $(u,v)\neq(v,u)$ $(e\in E, \, u,v\in V)$.
\begin{definition}\label{deg}
 Let $\mathcal{G}$ be a metric graph.

  \noindent $(i)$  The {\em degree} $d_{v}\in\N$ of a vertex $v$ is the number of edges $e\ni v$. The {\em inbound degree} $d_{v}^+$ (respectively, {\em outbound degree} $d_{v}^-$) is the number of edges with $j(e)=v$ (respectively, with $i(e)=v$).  $\mathcal{G}$ is {\em locally finite} if $d_v<\infty$ for any $v\in V$.

\noindent $(ii)$  The set $\pa \mathcal{G}:=\{v\in V\;|\; d_v=1\}$ is called the {\em boundary} of the graph. The vertices $v\in V\setminus \pa\mathcal{G}$ are called {\em interior}.

  \noindent $(iii)$ The {\em star centered at} $v$,
\begin{equation}\label{sv2}
 \Sig_v \,:=\,\{x\in\mathcal{G}\;|\; d(x,v)< l_j, \, j=1,\dots,d_v\}\qquad (v\in V)\,.
\end{equation}
We also define the {\em inbound star} $\Sig_v^+$ and the {\em outbound star $\Sig_v^-$ centered at} $v$ as the union of the edges ending at $v$, respectively emanating from $v$:
\begin{equation}\label{sv1}
\Sig_v^\pm \,:=\, \{x\in\mathcal{G}\;|\; d(x,v)< l_j, \, j=1,\dots,d_v^\pm\}\,,
\end{equation}
\end{definition}

Clearly, there holds $d_{v}=d_{v}^++d_{v}^-$ and $\Sig_v =\Sig_v^+(v) \cup \Sig_v^-(v)$  $(v\in V)$. It is not restrictive to assume
$\Sig_v^\pm\neq\emptyset$ for $v\in V\setminus\pa\mathcal{G}$.

\begin{definition}\label{sig}
$(i)$  Let $\mathcal{G}$ be a metric graph, and let $u,v\in V$. A {\em path} connecting $u$ and $v$ is a set $\{x_1,\dots,x_n\}\subseteq \mathcal{G}$ $(n\in\N)$ such that  $x_1=u$, $x_n=v$ and for all $k=1,\dots, n -1$ there exists an edge $e_k$ such that $x_k, x_{k+1} \in \overline{\mathcal{G}}_{e_k}$\,. A path connecting $u$ and $v$ is {\em closed} if $u\equiv v$. A closed path is called a {\em cycle}, if it does not pass through the same vertex more than once.

\noindent $(ii)$ A metric graph $\mathcal{G}$ is {\em connected} if for any $u,v\in V$ there exists a path connecting $u$ and $v$. A connected graph without cycles is called a {\em tree}.
\end{definition}

 It is easily seen that a connected metric graph $\mathcal{G}$ can be endowed with a structure of metric measure space. In fact, every two points $x,y\in\mathcal{G}$ can be regarded as vertices of a path $P$ connecting them  (possibly adding them to $V$). The length of the path is the sum of its $n$ edges $e_k$, i.e., $l(P):= \sum_{k=1}^n l_{e_k}\,,$ and the distance $d=d(x,y)$ between $x$ and $y$ is
$$
d(x,y)\,:=\,\inf\,\{l(P)\,|\, \text{$P$ connects $x$ and $y$}\}\,.
$$
Therefore $\mathcal{G}$ is a metric space, thus a topological space endowed with the metric topology.

Let $\mathcal{B}=\mathcal{B}(\mathcal{G})$ be the  Borel $\s$-algebra on $\mathcal{G}$. A Radon measure $\m:\mathcal{B}\mapsto [0,\infty]$ is induced on $\mathcal{G}$ by the Lebesgue measure $\l$ on each interval $I_e$, namely
\begin{equation}\label{demu}
\m(G)\,:=\, \sum_{e\in E} \l(I_e\cap G) \quad\text{for any $G\in\mathcal{B}$}\,.
\end{equation}

Let $\mathcal{G}$ be a metric graph. Every function $f:\mathcal{G}\mapsto\R$ canonically induces a countable family $\{f_e\}$, $f_e:I_e\mapsto\R$ $(e\in E)$. This is expressed as $f=\bigoplus_{e\in E}f_e$\,. We set $f^{(h)}:=\bigoplus_{e\in E}f_e^{(h)}$ $(h\in\N)$ if the derivative $f_e^{(h)}\equiv\frac{d^h f_e}{dx^h}$ exists in $I_e$ for all $e\in E$ (we shall also denote $f^{(0)}\equiv f$, $f^{(1)}\equiv f'$ and $f^{(2)}\equiv f''$). Let $C(\mathcal{G})$ denote the space of continuous functions on the metric space $(\mathcal{G},d)$. We set
$$
C^k(\mathcal{G})\,:=\,\{f\in C(\mathcal{G})\;|\; f_e\in C^k(I_e) \;\forall e\in E, \;\; f^{(h)}\in C(\mathcal{G})\;\forall \;h=1,\dots,k\} \qquad(k\in\N)\,,
$$
and $C^0(\mathcal{G})\equiv C(\mathcal{G})$\,. We also denote by $C^k_c(\mathcal{G})$ the subspace of functions in $C^k(\mathcal{G})$ with compact support, by $C_0(\mathcal{G})$ the closure of $C_c(\mathcal{G})\equiv C_c^0(\mathcal{G})$ with respect to the norm $\|\cdot\|_\infty$, and we set $C^\infty_c(\mathcal{G}):=\bigcap_{k\in\N}C^k_c(\mathcal{G})$\,.

In view of \eqref{demu}, for any measurable $f:\mathcal{G}\mapsto\R$ we set
$$
\int_\mathcal{G} f\,d\m\,:=\, \sum_{e\in E} \int_0^{l_e} f_e\,dx\,, \qquad \int_G f\,d\m\,:=\, \int_\mathcal{G} f\chi_G\,d\m
\quad\text{for any $G\in\mathcal{B}(\mathcal{G})$\,,}
$$
where $\chi_G$ denotes the characteristic function of the set $G$ and the usual notation $dx\equiv d\l$ is used. Accordingly, for any $p\in[1,\infty]$ the Lebesgue spaces $L^p(\mathcal{G})\equiv L^p(\mathcal{G},\m)$ are immediately defined:
$$
L^p(\mathcal{G})\,:=\,\bigoplus_{e\in E}\,L^p((I_e),\l)
$$
with norm
$$
\|f\|_p\,:=\,\sum_{e\in E}\|f_e\|_p\, =\, \sum_{e\in E} \left(\int_0^{l_e} |f_e|^p\,dx\right)^{\frac1p} \;\; \text{if $p\in (1,\infty)$}\,, \quad \|f\|_\infty\,:=\,{\rm ess\, sup}_{x\in\mathcal{G}}\, |f(x)|\,.
$$

For any $p\in[1,\infty]$ and $m\in\N$ the Sobolev spaces are
$$
W^{m,p}(\mathcal{G},\m)\,:=\,\left\{f\in L^p(\mathcal{G},\m)\;|\; f^{(h)}\in L^p(\mathcal{G},\m)\cap C(\mathcal{G})\;\forall \;h=0,\dots,m-1, \; f^{(m)}\in L^p(\mathcal{G},\m)  \right\}
$$
with norm $\|f\|_{m,p}\,:=\,\sum_{h=0}^m\|f^{(h)}\|_p$\,. As usual, we set $W^{m,p}(\mathcal{G})\equiv W^{m,p}(\mathcal{G},\m)$, $W^{0,p}(\mathcal{G})\equiv L^p(\mathcal{G})$,  $H^m(\mathcal{G}):= W^{m,2}(\mathcal{G})$, and we denote by $W^{m,p}_0(\mathcal{G})$ (in particular, by
$H^m_0(\mathcal{G})$) the closure of $C_c^\infty(\mathcal{G})$ with respect to the norm $\|\cdot\|_{m,p}$ (respectively $\|\cdot\|_{m,2}$). If $\mathcal{G}$ is compact, the embedding of $W^{1,p}(\mathcal{G})$ in $C_0(\mathcal{G})$ is continuous for $p\in[1,\infty)$ and compact for $p\in(1,\infty)$ (see \cite[Corollary 2.3]{H}).

For any $V_D\subseteq V$ we define
\begin{equation}\label{dewey}
W^{m,p}(\mathcal{G},V_D)\,:=\,\left\{f\in W^{m,p}(\mathcal{G})\;|\; f^{(h)}(v)=0 \;\forall\,v\in V_D, \;h=0,\dots,m-1 \right\},
\end{equation}
and $H^m(\mathcal{G},V_D)\equiv W^{m,2}(\mathcal{G},V_D)$. Clearly, there holds $W^{m,p}(\mathcal{G},\emptyset)=W^{m,p}(\mathcal{G})$,  and $W^{m,p}(\mathcal{G},\pa \mathcal{G})=W^{m,p}_0(\mathcal{G})$.   In particular, there holds
$$
H^1(\mathcal{G},V_D)\,=\,\left\{f\in H^1(\mathcal{G})\;|\; f(v)=0 \;\forall\,v\in V_D \right\},
$$
thus
$$
H^1(\mathcal{G},\pa\mathcal{G})\,=\,\left\{f\in H^1(\mathcal{G})\;|\;  f(v)=0 \;\forall\,v\in \pa\mathcal{G} \right\}\,=\,H_0^1(\mathcal{G})\,.
$$

\subsubsection{Laplacian on metric graphs}\label{lamegra}

Using the above functional framework a metric graph $\mathcal{G}$ can be endowed with a Dirichlet form, which makes it a  metric Dirichlet graph. The most obvious choice is the energy form
\begin{equation} \label{nufo}
\mathfrak{a}(f,g)\,:=\, \int_\mathcal{G} f'g'\,d\m\,,
\qquad\text{$f,g\in\mathfrak{D}(\mathfrak{a})\,:=\,H^1(\mathcal{G})$}
\end{equation}
(this choice is feasible, since the Hilbert space $H^1(\mathcal{G})$ is densely defined and continuously embedded in $L^2(\mathcal{G})$). The Laplacian $\D$ associated with $\mathfrak{a}$ is the operator with domain
\begin{subequations} \label{lave}
\begin{equation}\label{lave1}
D(\D)\,:=\, \left\{f \in H^1(\mathcal{G})\;|\;
f_e\in H^2(I_e) \,\forall e\in E, \; \sum_{e\in E}\int_0^{l_e}|f_e''|^2\,dx<\infty , \; \sum_{e\ni v}\frac{df_e}{d\n}(v)=0\,\forall\;v\in V\right\}
\end{equation}
(here $\frac{df_e}{d\n}(v)$ denotes the outer derivative of $f_e$ at the vertex $v$), which acts on $D(\D)$ in the natural way:
\begin{equation}
(\D f)_e\,:=\, f_e'' \quad \textrm{for any}\; f\in D(\D) \qquad(e\in E)\,.
\end{equation}
\end{subequations}

Since $H^1(\mathcal{G})\subseteq C(\mathcal{G})$, every $f\in D(\D)$ is continuous at any vertex $v\in \mathcal{G}\setminus  \pa\mathcal{G}$. The condition $\sum_{e\ni v}\frac{df_e}{d\n}(v)=0$ is the Kirchhoff transmission condition if $v\in \mathcal{G}\setminus \pa\mathcal{G}$, and the Neumann homogeneous boundary condition if $v\in \pa\mathcal{G}$\,. Hence the Laplacian defined in \eqref{lave} is called {\em Neumann Laplacian} and denoted by $\Delta_N$, if clarity so requires.

If $V_D\subseteq V$, $V_D\neq\emptyset$, a different notion of Laplacian is obtained considering the Dirichlet form
\begin{equation} \label{difo}
\mathfrak{a}_{_{V_D}}(f,g)\,:=\,\int_\mathcal{G} f'g'\,d\m\,,
\qquad\text{$f,g\in\mathfrak{D}(\mathfrak{a}_{_{V_D}})\,:=\,H^1(\mathcal{G},V_D)$.}
\end{equation}
The Laplacian $\D_{_{V_D}}$ associated with $\mathfrak{a}_{_{V_D}}$ is the operator defined as follows:
\begin{subequations} \label{lavedi}
\begin{equation}
(\D_{_{V_D}} f)_e\,:=\, f_e'' \quad \textrm{for any}\; f\in D(\D_{_{V_D}}) \quad(e\in E)\,,
\end{equation}
where
\begin{equation}
D(\D_{_{V_D}}):= \left\{f \in H^1(\mathcal{G},V_D)\,|\, f_e\in H^2(I_e) \,\forall e\in E,\,\sum_{e\in E}\int_0^{l_e}|f_e''|^2\,dx<\infty , \, \sum_{e\ni v}\frac{df_e}{d\n}(v)=0\,\forall\;v\in V\setminus V_D\right\}.
\end{equation}
\end{subequations}

More explicitly, the vertex conditions satisfied by all $f\in D(\D_{_{V_D}})\subseteq C(\mathcal{G})$ are:
$$
\left\{
\begin{array}{ll}
f(v)=0 &\quad\forall v\in V_D\,, \medskip \\
\sum_{e\ni v}\frac{df_e}{d\n}(v)=0 &\quad\forall v\in V\setminus V_D\ \, .
\end{array}
\right.
$$
Therefore, it is natural to call {\em Dirichlet Laplacian} the Laplacian $\D_D\equiv \D_{ \pa\mathcal{G}}$ (observe that the Neumann Laplacian satisfies $\D_N\equiv \D_\emptyset$). In particular, for a metric tree $\mathcal{T}$ the vertex conditions satisfied by any $f\in D(\D_D)\subseteq C(\mathcal{T})$ explicitly read:
$$
\left\{
\begin{array}{ll}
f(O)=0\,, & \medskip \\
\sum_{e\ni v}\frac{df_e}{d\n}(v)=0 &\quad\forall v\in \mathcal{T}\setminus O \, .
\end{array}
\right.
$$

\smallskip

Any realization of the Laplacian on a metric graph $\mathcal{G}$ ia a self-adjoint  nonpositive operator, thus the spectrum $\s(-\D)$ is contained in $\overline{\R}_+$.
We shall denote by $\{e^{\D t}\}_{t\ge0}$ the analytic semigroup generated by $\D$ in $L^2(\mathcal{G})$, which exists by classical results (e.g., see \cite{Y}). As in \cite[Theorems 1.4.1-1.4.2]{D}, the following holds:
\begin{prop}\label{davies}
The semigroup $\{e^{\D t}\}_{t\ge0}$ may be extended to a positivity preserving contraction semigroup $\{T_p(t)\}_{t\ge0}$ on $L^p(\mathcal{G})\equiv L^p(\mathcal{G},\m)$ for all $p \in[1,\infty]$. Moreover,  there holds
\begin{equation}\label{consi}
T_p(t)u= T_q(t)u \quad\text{for all $p,q \in[1,\infty]$ and $u \in L^p(\mathcal{G})\cap L^q(\mathcal{G})$} \,.
\end{equation}
The semigroup $\{T_p(t)\}_{t\ge0}$ is strongly continuous if $p \in[1,\infty)$ and holomorphic if $p \in(1,\infty)$.
\end{prop}
\noindent The notation $e^{\D t}\equiv T_p(t)$ $(p \in[1,\infty])$ will be used.

If for some smooth function $K:\mathcal{G}\times \mathcal{G}\times\overline{\R}_+\mapsto\overline{\R}_+$ there holds
\begin{equation}\label{edt}
e^{\D t}u\,=\int_\mathcal{G} K(\cdot,y,t)\,u(y)\,d\m(y) \qquad (u \in \bigcup_{p=1}^\infty L^p(\mathcal{G})),
\end{equation}
we say that the {\em heat kernel $K$ exists on} $\mathcal{G}$. Fundamental properties of the kernel are:

\noindent $(a)$ $K$ is symmetric, $K(x,y,t)=K(y,x,t)$, and smooth;

\noindent $(b)$ for any fixed $y\in \mathcal{G}$ $K(\cdot,y,\cdot)$ satisfies the heat equation in the classical sense,
\begin{equation}\label{caleq}
K_t(x,y,t)\,=\,\D_x K(x,y,t)\,;
\end{equation}

\noindent $(c)$ for any $x,y,z \in \mathcal{G}$ and $s,t \in \overline{\R}_+$
\begin{equation}\label{sgp}
\int_\mathcal{G} K(x,y,t)\,K(y,z,s)\,d\m(y)\,=\,K(x,z,s+t)\,;
\end{equation}

\noindent $(d)$ for any $x\in \mathcal{G}$ there holds
\begin{equation}\label{calid}
\lim_{t\to0^+} \int_\mathcal{G} K(x,y,t)\,u(y)\,d\m(y)\,=\, u(x)\quad\text{for all $u\in C(\mathcal{G})\cap L^\infty(\mathcal{G})$} \,.
\end{equation}

The existence of the heat kernel on metric graphs is ensured by \cite[Proposition 3.3]{KLVW}.
\smallskip

If $\mathcal{G}$ is compact, due to the compact embedding of $H^1(\mathcal{G})$ in $L^2(\mathcal{G})$, $\s(-\D)$ consists of eigenvalues $\l_j\ge0$ $(j\in \{0\}\cup\N)$ of finite multiplicity. The first eigenvalue $\l_0$ is 0 for $-\D_N$ and strictly positive for $-\D_D$\,. In both cases there exists an orthonormal basis $\{\varphi_j\}\subseteq L^2(\mathcal{G})$ of eigenfunctions such that
\begin{subequations}\label{spth}
\begin{equation}
D(-\D )\,=\, \left\{f\in L^2(\mathcal{G})\; \Big|\; \sum_{j=0}^\infty \l_j^2 \left(\int_\mathcal{G} f \varphi_j\,d\m\right)^2<\,\infty \right\},
\end{equation}
\begin{equation}
-\D f\,=\, \sum_{j=0}^\infty \l_j\left(\int_\mathcal{G} f \varphi_j\,d\m\right)\varphi_j \qquad (f\in D(-\D ))\,.
\end{equation}
\end{subequations}
By \eqref{spth} and the spectral theorem there holds
\begin{equation}\label{reker}
K(x,y,t)\,=\, \sum_{j=0}^\infty e^{-\l_j t} \varphi_j(x)\varphi_j(y)\,,
\end{equation}
where $K$ denotes the heat kernel either of $\D_N$ or of $\D_D$, and the series is uniformly convergent on $\mathcal{G}\times \mathcal{G}\times[\t,\infty)$ for any $\t>0$. Observe that for $\D_N$ equality \eqref{reker} reads
\begin{equation}\label{rekern}
K(x,y,t)\,=\, \frac{1}{\m(\mathcal{G})} \,+ \,\sum_{j=1}^\infty e^{-\l_j t} \varphi_j(x)\varphi_j(y)\,.
\end{equation}
We refer the reader to \cite[Proposition 2.6]{PT} for the proof of the following result.
\begin{prop}\label{cokerdi}
Let $\mathcal{G}$ be a metric graph, and let $\O_k\subseteq\mathcal{G}$ be an increasing family of subsets such that $\O_k\subset\subset \O_{k+1}$ $(k\in\N)$, $\bigcup_{k=1}^\infty \O_k=\mathcal{G}$. Let $\D$ be either the Dirichlet or the Neumann Laplacian on $\mathcal{G}$. Let $\D_k$ be either the Dirichlet or the Neumann Laplacian on $\O_k$, endowed with homogeneous Dirichlet boundary conditions on $\pa\O_k\setminus V$.

\smallskip

\noindent $(i)$ Let $K$, $K_k$ be the  attendant heat kernel of $\D$, respectively $\D_k$ $(k\in\N)$. Then there holds $K_k \le K_{k+1}\le  K$ in $\O_k$ for any $k\in\N$, and for any $u\in L^2(\mathcal{G})$
$$
\lim_{k\to\infty}e^{\D_k t}\left(u\chi_{_{\O_k}}\right)= e^{\D t}u\quad\text{ a.e. in $\mathcal{G}$ \qquad$(t>0)$\,.}
$$

\noindent $(ii)$ Let  $-\l_{0,k}<0$ be the first eigenvalue of $\D_k$ $(k\in\N)$, and $E_0:=\min\s(-\D)\ge0$\,. Then there holds $-\l_{0,k}\le -\l_{0,k+1}\le -E_0$ for any $k\in\N$, and $\lim_{k\to\infty}\l_{0,k}\,=\,E_0\,$.
\end{prop}


\subsection{Metric trees}\label{metre}
In this subsection we first recall some general concepts concerning metric trees. Further, we describe a characterization of the Laplacian restricted to the subspace of symmetric functions,  which plays an important role to prove spectral properties of the Laplacian itself.

\subsubsection{General facts}\label{gefatre}

\begin{definition}\label{detree}
\noindent $(i)$ A {\em metric tree} is a connected metric graph $\mathcal{T}$ without cycles, with a boundary vertex $O$ singled out. $O$ is called the {\em root} of $\mathcal{T}$. The quantity $h(\mathcal{T}):=\sup_{x\in\mathcal{T}}\r(x)$\,, where $\r(x):=d(x,O)$ $(x\in \mathcal{T})$, is called the {\em height} of $\mathcal{T}$.

\noindent $(ii)$ The {\em generation} ${\rm gen}(v)$ of a vertex $v \in \mathcal{T}$ is the number of vertices which lie on the unique path connecting $v$ with $O$ (including the starting point but excluding the end point). The {\em branching number} $b(v)$ of $v \in \mathcal{T}$ is the number of edges emanating from $v$.

\noindent $(iii)$ A tree $\mathcal{T}$ is {\em regular}, if all vertices of the same generation have equal branching numbers, and all edges emanating from them are of the same length. Clearly,

\noindent $(a)$ in a regular tree all vertices of the same generation have the same distance from the root;

\noindent $(b)$ a regular tree is uniquely determined by two {\em generating sequences} $\{b_n\}\equiv\{b_n(\mathcal{T})\}\subseteq\{0\}\cup\N$,
$\{\r_n\}\equiv\{\r_n(\mathcal{T})\}\subseteq\R_+$, with $\{\r_n\}$ increasing and $\r_0=0$, such that
$$
n={\rm gen}(v)\,, \quad b_n=b(v)\,,\quad \r_n=\r(v) \quad\text{for any vertex $v$ of $\mathcal{T}$}
$$
(namely, $b_n$ is the common branching number of all vertices of generation $n$).

\noindent $(iv)$ A regular tree $\mathcal{T}$ is {\em homogeneous}, if the branching number $b(v)$ is the same for all $v \in \mathcal{T}\setminus O$, and all edges are of the same length.
\end{definition}

We always suppose that $b_0 =1$, $b_n\ge2$ for any $n\in\N$, ${\rm gen}(v)<\infty$ for all $v \in \mathcal{T}$, and $h(\mathcal{T})=\infty$\,. In particular, there holds $\pa\mathcal{T}=O$ (i.e., the boundary of a metric tree consists of its root). Observe that $\bigcup_{n\in\N}(\r_{n-1},\r_n)=(0,h(\mathcal{T}))=\R_+$\,.
If all edges of a regular tree have the same length $r>0$, for all $n\in\N$ there holds $\r_n- \r_{n-1}=r$, thus $\r_n=nr$.

If the tree $\mathcal{T}$ is regular, its {\em branching function} $\b:[0,h(\mathcal{T}))\mapsto \N$ is defined as follows:
\begin{equation}\label{bratree}
\b(\r):={\rm card}\{x\in\mathcal{T}\,|\, \r(x)=\r \}\,=\, \left\{
\begin{array}{ll}
1  &\text{if $\r=0$\,,} \\
\sum_{n=1}^\infty ( b_0b_1\dots b_{n-1})\chi_{(\r_{n-1},\r_n]}(\r)&\text{if $\r>0$}\,.
\end{array}
\right.
\end{equation}
Clearly, $\b$ is an increasing step function continuous from the left, with jump $\b(\r_n^+)-\b(\r_n^-)=b_0b_1\dots b_{n-1}(b_n-1)$ at any point $\r=\r_n$ $(n\in\N)$. Observe that $\b(\r_n)\ge 2^{n-1}$, whereas $\b(\r_n)=b^{n-1}$ for a homogeneous tree ($n\in\N$).

Therefore, if the length of all edges of a regular tree exceeds some $r_0>0$,
the volume $V(O,R):=\m(B(O,R))$ of the ball $B(O,R):=\{x\in\mathcal{T}\,|\,\r(x)<R\}$ satisfies the inequalities
\begin{equation}\label{voltree}
V(O,R) \,\ge\, V(O,[R])\,\ge\, r_0\sum_{n=1}^{[R]}\b(\r_n) \,\ge\,  r_0\, 2^{[R]-1} \,,
\end{equation}
where $[R]$ denotes the largest integer not exceeding $R$. Hence in this case the volume growth of $\mathcal{T}$ is exponential.


\subsubsection{Symmetric functions}\label{ss32}

As pointed out in the Introduction, the special structure of regular trees provides much information about the  spectral properties of the Laplacian on them. In this subsection we give an outline of general results in this direction, referring the reader to \cite {NS1,NS2,So1,So2} for the proofs. Let us first state the following definition.
\begin{definition}
Let $\mathcal{T}$ be a regular tree. A function $f:\mathcal{T}\mapsto [0,1]$ is called {\em symmetric} if
$$
x, y \in \mathcal{T},\;\,   \r(x)=\r(y) \quad  \Rightarrow \quad f(x)=f(y)\,.
$$
\end{definition}
Clearly, $f$ is symmetric if and only if there exists $\tilde f:[0,h(\mathcal{T}))\mapsto [0,1]$ such that
\begin{equation}\label{ftf}
f(x)=\tilde f(\r) \quad\text{for any $x\in \mathcal{T}$ with $\r(x)=\r$}\,,
\end{equation}
namely such that
\begin{equation}\label{ftf1}
f\,=\,\tilde f\circ\r \,.
\end{equation}
\begin{remark}\label{marko0}
Since the map $x\mapsto\r(x)$ $(x\in\mathcal{T})$ is continuous, by \eqref{ftf} for any $f\in C(\mathcal{T})$ there holds $\tilde f\in C(0,h(\mathcal{T}))$.
\end{remark}
The following result is immediately proven (see \cite{NS2}).
 \begin{lem}
 Let $\mathcal{T}$ be a regular tree. Let  $f:\mathcal{T}\mapsto [0,1]$ be symmetric, and let $\tilde f:[0,h(\mathcal{T}))\mapsto [0,1]$
be as in \eqref{ftf}. Then
\begin{equation}\label{leftf}
\sum_{x\in\mathcal{T}\!\!,\,\r(x)=\r} f(x)\,=\, \tilde f(\r)\,\b(\r) \quad\text{for a.e. $x\in\mathcal{T}\,.$}
\end{equation}\,
 \end{lem}

Set
\begin{equation}\label{dl2b}
L^2(\R_+;\b)\,:=\, \left\{\tilde f:\R_+\mapsto\R\;|\; \tilde f \,\textrm{\,measurable}, \int_{\R_+}|\tilde f|^2\,\b\,d\r<\infty \right\}\,,
\end{equation}
and denote by $H^1(\R_+;\b)$ the space of functions $\tilde f\in L^2(\R_+;\b)$ with distributional derivative $\tilde f'\in L^2(\R_+;\b)$\,. Both $L^2(\R_+;\b)$ and $H^1(\R_+;\b)$ are Hilbert spaces with norms
$$
\|\tilde f\|_{L^2\!,\, \b}\,:=\,\left ( \int_{\R_+} |\tilde f|^2\,\b\,d\r \right)^{\frac12}\,,
$$
respectively
$$
\|\tilde f\|_{H^1\!,\, \b}\,:=\,\left ( \int_{\R_+} \left [\,|(\tilde f)'|^2+|\tilde f|^2\right]\,\b\,d\r \right)^{\frac12}\,.
$$

Let $\mathfrak{S}(\mathcal{T})$ denote the subspace of symmetric functions on $\mathcal{T}$ (see Definition \ref{detree}-$(i)$). Set
$$
\mathfrak{L}^2(\mathcal{T})\,:=\, L^2(\mathcal{T}) \cap\mathfrak{S}(\mathcal{T})\,,  \qquad \mathfrak{H}^1(\mathcal{T})\,:=\, H^1(\mathcal{T}) \cap\mathfrak{S}(\mathcal{T})\,.
$$
In view of \eqref{leftf}, there holds
\begin{equation}\label{l2ftf}
\|f\|_2^2\,\equiv\, \|f\|_{L^2(\mathcal{T})}^2\,=\, \int_\mathcal{T} |f|^2\,d\m\,=\, \int_{\R_+}|\tilde f|^2\,\b\,d\r\,=\,\|\tilde f\|_{L^2\!,\, \b}^2\,<\,\infty \quad\text{for all $f\in\mathfrak{L}^2(\mathcal{T})$}\,,
\end{equation}
hence $\mathfrak{L}^2(\mathcal{T})$ is isometrically isomorphic to the weighted space $L^2(\R_+;\b)$. If $f\in \mathfrak{H}^1(\mathcal{T})$, thus $f'\in \mathfrak{L}^2(\mathcal{T})$, it is easily checked that $\tilde{(f')}=(\tilde f)'$ a.e. in $\R_+$\,, thus by \eqref{l2ftf} there holds
\begin{equation}\label{h1ftf}
\int_\mathcal{T} |f'|^2\,d\m=\, \int_{\R_+}|\tilde{(f')}|^2\,\b\,d\r\,=\, \int_{\R_+}|(\tilde f)'|^2\,\b\,d\r\,<\,\infty \quad\text{for all $f\in\mathfrak{H}^1(\mathcal{T})$}\,.
\end{equation}
By \eqref{l2ftf}-\eqref{h1ftf}, for all $f\in\mathfrak{H}^1(\mathcal{T})$
$$
\|f\|_{1,2}^2\equiv \|f\|_{H^1(\mathcal{T})}^2 \,=\, \int_\mathcal{T} \left[\,[|f'|^2+|f|^2\right] \,d\m\,=\,  \int_{\R_+} \left [\,|(\tilde f)'|^2+|\tilde f|^2\right]\,\b\,d\r \,=\,\|\tilde f\|_{H^1\!,\, \b}^2\,,
$$
thus $\mathfrak{H}^1(\mathcal{T})$ is isometrically isomorphic to the weighted space $H^1(\R_+;\b)$\,.
\begin{remark}\label{marko}
Since $\mathfrak{H}^1(\mathcal{T})\subseteq H^1(\mathcal{T})\subseteq C(\mathcal{T})$, by Remark \ref{marko0} 
for any $f\in \mathfrak{H}^1(\mathcal{T})$ there holds $\tilde f \in C(\R_+)$\,.
\end{remark}
Let $\mathfrak{a}$ be the energy form which defines the Neumann Laplacian $\D$ on $\mathcal{T}$ (see \eqref{nufo}):
\begin{equation*}
\mathfrak{a}(f,g)\,:=\, \int_\mathcal{T} f'g'\,d\m\,,
\qquad \mathfrak{D}(\mathfrak{a})\,:=\,H^1(\mathcal{T})\,.
\end{equation*}
If $f,g\in\mathfrak{H}^1(\mathcal{T})= H^1(\mathcal{T}) \cap\mathfrak{S}(\mathcal{T})$, the same form $\mathfrak{a}(f,g)$ defines the restriction $\D\big|_{\mathfrak{S}(\mathcal{T})}$, and there holds
\begin{subequations}\label{mnz}
\begin{equation}
\mathfrak{a}(f,g)\,:=\, \int_\mathcal{T} f'g'\,d\m\,=\,\int_{\R_+} (\tilde f)'(\tilde g)'\b\,d\r\,=:\,\tilde{\mathfrak{a}}(\tilde f,\tilde g)\,,
\end{equation}
the domain of the form $\tilde{\mathfrak{a}}$ being
\begin{equation}
D(\tilde{\mathfrak{a}})\,:=\,H^1(\R_+;\b)\,.
\end{equation}
\end{subequations}

Now consider the unitary transformation
\begin{equation}\label{truni}
U: L^2(\R_+;\b) \mapsto L^2(\R_+)\,,\qquad \tilde f \mapsto F\equiv U\tilde f:= \sqrt{\b}\,\tilde f\,.
\end{equation}
Since the branching function $\b$ is constant in each interval $I_n:=(\r_{n-1},\r_n)$ $(n\in\N)$, where $\{\r_n\}$ is a generating sequence of $\mathcal{T}$, there holds $F'(\r)= \sqrt{\b}(\r)(\tilde f)'(\r)$ for a.e. $\r\in\R_+$\,. Therefore, the form $\tilde{\mathfrak{a}}$ defined in \eqref{mnz} corresponds to the form
\begin{subequations}\label{mnbbis}
\begin{equation}
\mathfrak{a}_0(F,G)\,:=\, \int_{\R_+} F'G'\,d\r\,,
\end{equation}
with domain
\begin{equation}\label{emme}
D(\mathfrak{a}_0)\,:=\,\left\{F \in L^2(\R_+) \,\Big|\,F=\sum_{n=1}^\infty F_n\chi_{_{I_n}},\,F_n\in H^1(I_n),\sum_{n=1}^\infty\int_{I_n} |F_n'|^2\,d\r<\infty, F(\r_n^+)= \sqrt{b_n}\,F(\r_n^-)\,\forall n\in\N  \right \}
\end{equation}
\end{subequations}
(here $\chi_{_{I_n}}$ denotes the characteristic function of the interval $I_n$).
\begin{remark}\label{recon1}
The matching condition in \eqref{emme} stems from the continuity of $\tilde f$ in $\R_+$, since
$$
F(\r_n^+):= \sqrt{\b(\r_n^+)}\,\tilde f(\r_n)= \sqrt{b_n}\sqrt{b_1\dots b_{n-1}}\,\tilde f(\r_n) =\sqrt{b_n}\sqrt{\b(\r_n^-)}\,\tilde f(\r_n)= \sqrt{b_n}\,F(\r_n^-)\,.
$$
\end{remark}

The operator $A$ associated to $\mathfrak{a}_0$ by the Friedrichs' construction is easily described (e.g., see \cite[Theorem VI.2.1]{K}). Let $F\in D(\mathfrak{a}_0)$. Observe preliminarily that, if there exists $H\in L^2(\R_+)$ such that
\begin{subequations}\label{fried}
\begin{equation}\label{friedb}
\mathfrak{a}_0(F,G)\,=\, \int_{\R_+} GH\,d\r\quad\textrm{for all $G$ in a core of  $D(\mathfrak{a}_0)$}\,,
\end{equation}
then there holds $F\in D(A)$ and $AF=H$. Conversely, for any $F\in D(A)$
\begin{equation}\label{frieda}
\mathfrak{a}_0(F,G)\,=\, \int_{\R_+} (A F)G\,d\r\quad\textrm{for all $G\in D(\mathfrak{a}_0)$}\,.
\end{equation}
\end{subequations}

\begin{prop}\label{colf}
\begin{subequations}\label{frado}
Let $\mathcal{T}$ be a regular tree, and let $A$ be the operator associated to the form $\mathfrak{a}_0$ defined in \eqref{mnbbis}. Then the domain of $A$ is the set of functions $F$ with the following properties:
\begin{equation}\label{frado1}
F\in L^2(\R_+)\,, \qquad F_n\in H^2(I_n)\quad\text{for all $n\in\N$}\,,\qquad \sum_{n=1}^\infty\int_{I_n}\big(|F_n''|^2 +|F_n'|^2\big)\,d\r\,<\,\infty\,,
\end{equation}
\begin{equation}\label{frado2}
F(\r_n^+)\,=\, \sqrt{b_n}\,F(\r_n^-)\quad\text{for all $n\in\N$}\,,
\end{equation}
\begin{equation}\label{frado3}
F'(0^+)\,=\,0\,,\qquad F'(\r_n^-)\,=\, \sqrt{b_n}\,F'(\r_n^+)\quad\text{for all $n\in\N$}\,,
\end{equation}
\end{subequations}
with $\{b_n\}$ generating sequence of $\mathcal{T}$. The operator $A$ is self-adjoint and nonnegative, and acts on
$D(A)$ as follows:
\begin{equation}\label{supro3}
AF\,=\, -F''\quad\text{for any $F\in D(A)$\,.}
\end{equation}
\end{prop}

\begin{proof}
Let $G\in C_c([0, \r_1))\cap C^1(\bar I_1) \subseteq D(\mathfrak{a}_0)$. Then there holds
\begin{subequations}\label{pinte}
\begin{equation}
\mathfrak{a}_0(F,G)\,=\, \int_0^{\r_1} F'G'\,d\r
\,=\, -\int_0^{\r_1} F''G\,d\r \,-\, F'(0^+)\b(0)\,.
\end{equation}
Similarly, if $G\in \big[C_c((\r_{n-1},\r_n])\cap C^1(\bar I_n)\big]\cup \big[C_c([\r_n, \r_{n+1}))\cap C^1(\bar I_{n+1})\big] \subseteq D(\mathfrak{a}_0)$ $(n\in\N)$, there holds
\begin{equation}
\mathfrak{a}_0(F,G)\,=\, \int_{\r_{n-1}}^{\r_{n+1}} F'G'\,d\r
\,=\, -\int_{\r_{n-1}}^{\r_{n+1}}F''G\,d\r \,+\, \left[F'(\r_n^-)- \sqrt{b_n}\,F'(\r_n^+)\right]G(\r_n^-)
\end{equation}
\end{subequations}
(here use has been made of the matching conditions in \eqref{emme} for $G$; observe that the limits $F'(\r_n^\pm):=\lim_{\r\to\r_n^\pm}F'(\r)$ exist and are finite, since $F_n\in H^2(I_n)\subseteq C^1(\bar{I}_n)$  by embedding results).

Let $F=\sum_{n=1}^\infty F_n\chi_{_{I_n}}\in D(\mathfrak{a}_0)$ satisfy
\begin{subequations}\label{supro}
\begin{equation}\label{supro1}
F_n\in H^2(I_n)\quad\text{for all $n\in\N$}\,,\quad \int_{\R_+}|F''|^2\,d\r\,=\,\sum_{n=1}^\infty\int_{I_n}|F_n''|^2\,d\r\,<\,\infty\,.
\end{equation}
Then by proper choices of $G$, by \eqref{friedb} and  \eqref{pinte} we obtain that $F\in D(A)$ and $AF\,=\, -F''$, if there holds
\begin{equation}\label{supro2}
F'(0^+)\,=\,0\,,\quad F'(\r_n^-)\,=\, \sqrt{b_n}\,F'(\r_n^+)\quad\text{for all $n\in\N$}\,.
\end{equation}
\end{subequations}
Therefore, every function $F$ with properties \eqref{supro1}-\eqref{supro2} belongs to $D(A)$. Conversely, by \eqref{frieda} and the very definition of distributional derivative, for any $F\in D(A)$ there holds $AF\,=\, -F''\in L^2(\R_+)$, thus the above argument can be inverted to prove properties \eqref{supro1}-\eqref{supro2}.
Recalling that $D(A)\subseteq D(\mathfrak{a}_0)$, the conclusion follows.
\end{proof}

In view of the above discussion, the restriction of the (Neumann) Laplacian $\Delta$ to  the subspace $\mathfrak{S}(\mathcal{T})$ of symmetric functions is characterized as follows (see  \cite[Theorem 4.1]{NS2}, \cite[Lemma 3.4]{So2}).
\begin{prop}\label{clf}
Let $\mathcal{T}$ be a regular tree. Then the restriction $(-\Delta)\big|_\mathfrak{S(\mathcal{T})}$ is unitarily equivalent to the operator $A$ defined by the form \eqref{mnbbis}. The domain and the action of $A$ are described by \eqref{frado} and \eqref{supro3}, respectively.
\end{prop}
\begin{remark}\label{recon2}
In \eqref{frado3} the first equality corresponds to the Neumann boundary condition at the root $O$, whereas the jump condition at $\r_n$ stems from the Kirchhoff condition at any vertex $v$ with ${\rm gen}(v)=n$. In fact, in this case the Kirchhoff condition $\sum_{e\ni v}\frac{df_e}{d\n}(v)=0$ reads $f'(v^-)\, =\, b_n f'(v^+)$, namely
$$
(\tilde f)'(\r_n^-) \, =\, b_n (\tilde f)'(\r_n^+) \qquad \Longleftrightarrow \qquad  \frac{F'(\r_n^-)}{\sqrt{\b(\r_n^-)}}\, =\, b_n \, \frac{F'(\r_n^+)}{\sqrt{\b(\r_n^+)}}
$$
(see \eqref{truni}). Since $\b(\r_n^+)=b_n\b(\r_n^-)$, the claim follows.
\end{remark}



\end{document}